\documentclass[12pt]{article}
\usepackage[nosumlimits]{amsmath}
\usepackage{amssymb,amsthm,MnSymbol,url}

\def\MM#1{\boldsymbol{#1}}
\DeclareMathOperator{\diff}{d\!}

\usepackage{amscd}
\usepackage{natbib}
\usepackage{times}
\usepackage{amsfonts}
\usepackage{algorithmicx}
\usepackage[noend]{algpseudocode}
\usepackage{algorithm}

\usepackage[margin=2cm]{geometry}
\usepackage{graphicx}
\usepackage{color}
\usepackage{caption}
\usepackage{subcaption}
\usepackage{bm}
\usepackage{empheq}
\usepackage{authblk}

\usepackage[title]{appendix}

\begin{document}
\title{Conservation with moving meshes over orography}

\author[1$\dag$*]{Hiroe Yamazaki}
\author[1]{Hilary Weller}
\author[2]{Colin J. Cotter}
\author[3]{Philip A. Browne}

\affil[1]{Department of Meteorology, University of Reading, Reading, UK}
\affil[2]{Department of Mathematics, Imperial College London, London, UK}
\affil[3]{European Centre for Medium-Range Weather Forecasts, Reading, UK}
\affil[$\dag$]{Current affiliation: Department of Mathematics, Imperial College London, London, UK}
\affil[*]{Correspondence to: \texttt{h.yamazaki@imperial.ac.uk}}

\maketitle

\begin{abstract}
Adaptive meshes have the potential to improve the accuracy and efficiency of atmospheric modelling by increasing resolution where it is most needed. Mesh re-distribution, or r-adaptivity, adapts by moving the mesh without changing the connectivity. This avoids some of the challenges with h-adaptivity (adding and removing points): the solution does not need to be mapped between meshes, which can be expensive and introduces errors, and there are no load balancing problems on parallel computers. A long standing problem with both forms of adaptivity has been changes in volume of the domain as resolution changes at an uneven boundary. We propose a solution to exact local conservation and maintenance of uniform fields while the mesh changes volume as it moves over orography. This is solved by introducing a volume adjustment parameter which tracks the true cell volumes without using expensive conservative mapping.

A finite volume solution of the advection equation over orography on moving meshes is described and results are presented demonstrating improved accuracy for cost using moving meshes. Exact local conservation and maintenance of uniform fields is demonstrated and the corrected mesh volume is preserved.

We use optimal transport to generate meshes which are guaranteed not to tangle and are equidistributed with respect to a monitor function. This leads to a Monge-Amp\`{e}re equation which is solved with a Newton solver. The superiority of the Newton solver over other techniques is demonstrated in the appendix. However the Newton solver is only efficient if it is applied to the left hand side of the Monge-Amp\`{e}re equation with fixed point iterations for the right hand side.
\end{abstract}

\noindent \textbf{keywords:} moving meshes; mesh re-distribution; r-adaptivity; numerical weather prediction; orography

\section{Introduction}

Dynamic mesh adaptivity can be advantageous for the numerical solution of PDEs when numerical errors (or their impacts) are greater in some areas than others.
Numerical weather and climate predictions, for example, could be improved by locally varying the spatial resolution through time, tracking atmospheric phenomena such as weather fronts \cite[]{budd2013monge}.

$H$-adaptivity involves adding and removing computational points based on local resolution requirements \citep[e.g.][]{berger1984adaptive, skamarock1993adaptive,weller2009predicting}.
The connectivity of the mesh and the total number of computational points change.
Conversely, $r$-adaptivity, or mesh re-distribution, involves moving mesh vertices without changing the connectivity of the mesh.
It results in a deformed mesh keeping the number of computational points and topology the same \citep[e.g.][]{dietachmayer1992application, hirt1997arbitrary, kuhnlein2012modelling}.

$R$-adaptivity is an attractive form of adaptivity since the data structures associated with the connectivity do not change and therefore the load balancing remains constant on parallel computers.
$R$-adaptivity does not require mapping solutions between old and new meshes.
$R$-adaptivity can also lead to smoothly graded meshes, which are desirable in order to reduce wave reflections and other errors associated with rapid resolution changes \citep[][]{vichnevetsky1987wave, long2011numerical}.

\begin{figure}[t]
  \centering
  \includegraphics[height=45mm]{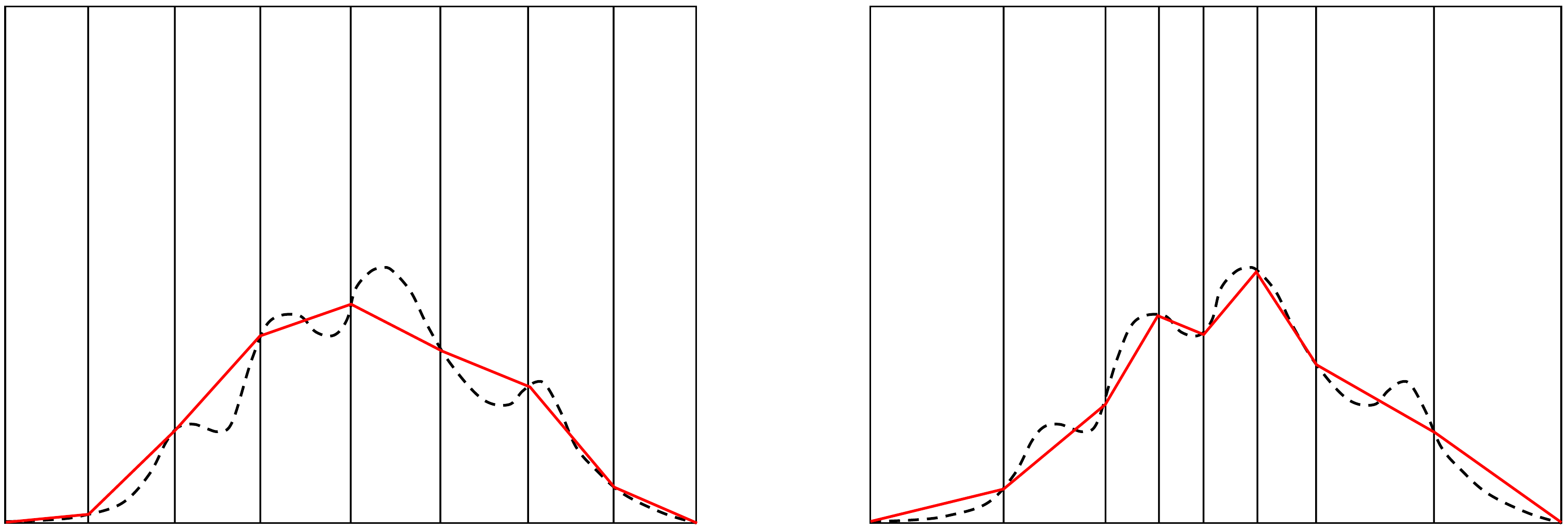} 
  \caption{Two-dimensional example of a change in shape of orography when the mesh distribution changes from (left) uniform to (right) focused in the middle. Red and dashed lines indicate the model orography and real orography, respectively.}
  \label{fig:topography}
\end{figure}

A disadvantage of $r$-adaptivity is, with a fixed number of points and fixed connectivity, it is not possible to achieve exactly the required resolution in each direction simultaneously, which can be achieved with $h$-adaptivity.
This causes an extra difficulty when introducing orography.
For example, we consider a vertical slice model with a horizontally moving mesh over orography (Figure \ref{fig:topography}).
When the topographic surface is approximated as piecewise linear (or higher order) splines based on vertex locations, the shape of orography inevitably changes as mesh vertices move over the orography.
Therefore the volume of the domain changes as the mesh moves, which results in unphysical compression or expansion of the model fluid.

With $h$-adaptivity, this issue can be avoided by evaluating all the metric terms from orography on the finest grid everywhere in the domain and average the result, in a consistent manner, to any coarser grid \citep[][]{guzik2015freestream}.
This approach is only suitable for models in which each fine grid cell only ever over-laps with one coarse grid cell.
Another way to resolve this issue is using conservative mapping to calculate the cell volumes over the original shape of orography.
For example, \cite{schwartz2015high} presented an algorithm to perform highly accurate mappings using sub-grid knowledge, which ensures the required accuracy with grid refinement.
Though this approach can be used with both $h$- and $r$-adaptivities, it is expensive to perform conservative mapping of the orography every time step.

Instead of tracking the shape of orography within each cell, we propose another solution which is to calculate the ``true'' cell volumes indirectly by solving a transport equation for the cell volume.
This is solved by introducing a volume adjustment parameter which tracks the change in cell volumes caused by the change in the shape of orography.
With this approach, the exact local conservation and maintenance of uniform fields on a moving mesh over orography is achieved without using expensive conservative mapping.

Section 2 provides the model description, including the finite-volume discretisation on a moving mesh, and the adjustment of the cell volumes as mesh moves over orography.
In section 3, we present the results of a three-dimensional tracer advection test with the use of the volume adjustment parameter.
Here we evaluate the model error using smooth and steep orography and demonstrate the importance of maintaining uniform fields on a moving mesh.
Finally, in section 4 we provide a summary and outlook.
We prove that the method for calculating adjusted cell volumes is bounded in appendix \ref{appx:bounded} and we describe the optimally transported mesh generation in appendix \ref{secn:MAsolution}.

\section{Model Description}\label{model_description}

\subsection{Finite volume discretisation on a moving mesh}
\label{discretisation}

\begin{figure}[h]
  \centering
  \includegraphics[height=60mm]{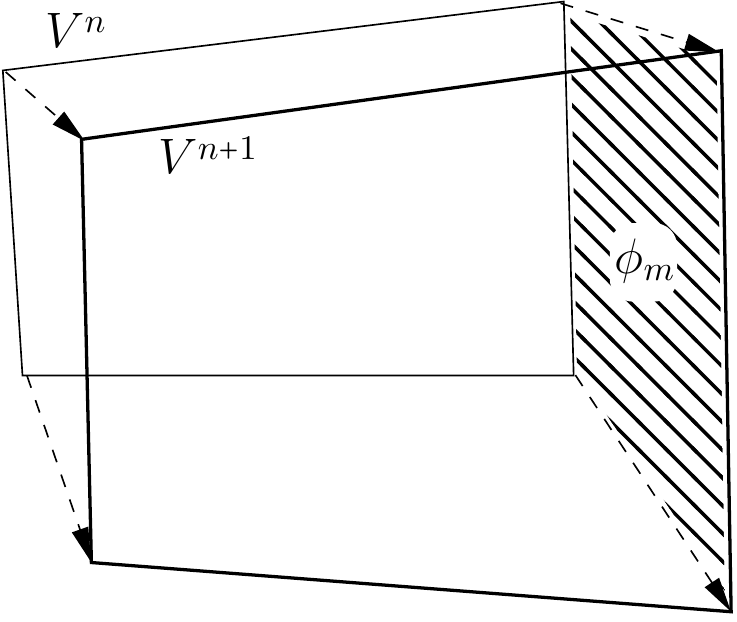} 
  \caption{An example showing how the mesh flux $\phi_m$ is calculated in the model. The solid line and the thick line show the boundaries of a cell at the time step $n$ and $n+1$, respectively. The hatched region shows the swept volume by the right face of the cell, which is calculated as the mesh flux at time step $n + 1/2$ corresponding to that face.}
  \label{fig:mesh_flux}
\end{figure}

We consider the three-dimensional advection equation in flux form:
\begin{eqnarray}
\frac{\partial \rho}{\partial t} + \nabla \cdot (\MM{u}\rho) &=& 0, \label{original_equation}
\end{eqnarray}
where $\MM{u} = (u, v, w)$ is a prescribed velocity field and $\rho$ is the tracer density.
To derive a finite-volume descretised equation, first we integrate the equation over a control volume $V$ and then apply Gauss' divergence theorem:
\begin{eqnarray}
\int_{V}\frac{\partial \rho}{\partial t} \diff V
+
\oint_{S} \rho \, \MM{u} \cdot \MM{n} \diff S
&=& 0,
\label{integral_form}
\end{eqnarray}
where $\MM{n}$ is the outward pointing unit normal vector on the boundary surface $S$ of the control volume $V$ so that $\rho \, \MM{u} \cdot \MM{n} \diff S$ is the flux of $\rho$ over the surface with area $\diff S$.
To extend \eqref{integral_form} to a moving mesh, we use the Reynolds transport theorem:
\begin{eqnarray}
\frac{\diff}{\diff t}\int_{V(t)} \rho \diff V = \int_{V(t)} \frac{\partial \rho}{\partial t} \diff V + \oint_{S(t)}\rho \, \MM{u_{s}} \cdot \MM{n} \diff S, \label{reynolds_theorem}
\end{eqnarray}
where $\MM{u_s}$ is the velocity of the boundary surface $S$.
Note that the control volume $V$ and the boundary surface $S$ are now time dependent.
The relationship between the volume $V$ and the velocity $\MM{u_s}$ is called the space conservation law \citep[][]{demirdvzic1988space}:
\begin{eqnarray}
\frac{\partial}{\partial t}\int_{V(t)} \diff V - \oint_{S(t)} \MM{u_s}\cdot\MM{n} \diff S = 0.
\label{space_conservation_law}
\end{eqnarray}
This means that the change in the control volume has to match the sum of the swept volumes of all its surfaces.
Combining the equations \eqref{integral_form} and \eqref{reynolds_theorem}, we have the integral form of the advection equation \eqref{original_equation} on a moving mesh:
\begin{eqnarray}
\frac{\diff}{\diff t}\int_{V(t)} \rho \diff V + \oint_{S(t)}\rho \, (\MM{u} - \MM{u_{s}}) \cdot \MM{n} \diff S = 0.
\label{moving_mesh_form}
\end{eqnarray}
The discretised form of the equation \eqref{moving_mesh_form} can be written:
\begin{eqnarray}
\frac{V^{n+1}\rho^{n+1}-V^{n}\rho^{n}}{\Delta t} + \sum\limits_{\rm faces}\rho_{f}(\phi - \phi_{m}) = 0,
\label{discretised_form}
\end{eqnarray}
where values inside the summation are at time step $n+1/2$, $\Delta t$ is the time step size, $\rho_{f}$ denotes the tracer density that is interpolated onto the cell faces, and $\phi = \MM{u}\cdot\MM{n}\diff S$ is the face flux defined on the cell faces.
The mesh flux $\phi_m = \MM{u_s}\cdot\MM{n} \diff S$ is calculated as the swept volume by the faces during each time step (Figure \ref{fig:mesh_flux}), which satisfies the following discretised form of the equation \eqref{space_conservation_law}:
\begin{eqnarray}
\frac{V^{n+1}-V^{n}}{\Delta t} - \sum\limits_{\rm faces}\phi_{m} = 0.
\label{mesh_flux}
\end{eqnarray}

\subsection{Automatic mesh motion}
\label{moving_mesh}

The optimally transported mesh generation procedure is described in appendix \ref{secn:MAsolution}. The Monge-Amp\`ere equation is solved to generate a mesh that is equidistributed with respect to a monitor function and guaranteed tangle free due to the optimal transport. A Newton method is described to solve the Monge-Amp\`ere equation. A monitor function is chosen so that the cell areas are a factor of 4 smaller in regions where the second derivatives of the tracer density is highest compared with the regions of lowest second derivatives. The mesh is moved every time step so fast convergence of the Newton solver is important, which is also demonstrated in appendix \ref{secn:MAsolution}.

\subsection{Volume adjustment as mesh moves over orography}
\label{colin_parameter}
When the cell volumes are calculated from vertex locations without tracking the variations in orography within each cell,
the shape of the model orography inevitably changes as the mesh vertices move (Figure \ref{fig:topography}).
This means that equation \eqref{mesh_flux} does not hold on a moving mesh over orography: we only consider the swept volumes in the horizontal, not the swept volume of the orography surface.
The mismatch between the change in the control volume and the sum of the swept volumes results in non-uniform model fields at the bottom surface.

\begin{figure}[t]
  \centering
  \includegraphics[height=45mm]{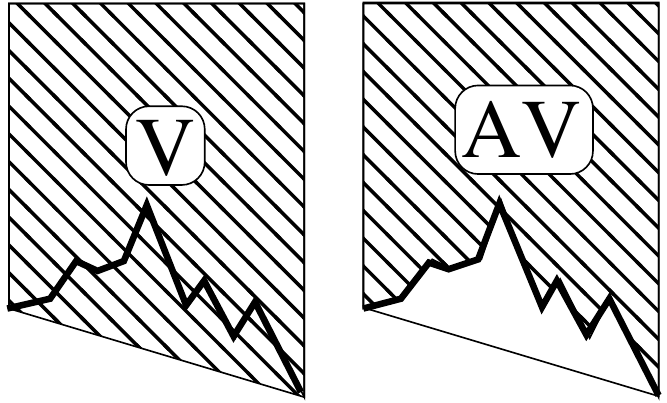} 
  \caption{Two-dimensional example of (left) the cell volume $V$ calculated from vertex locations, and (right) the true cell volume which is given by multiplying the volume adjustment parameter $A$ by $V$. The thick line represents the surface of orography. Hatched regions describe the volume of the cells.}
  \label{fig:AV}
\end{figure}

One solution to this problem would be to use a conservative mapping of old to new mesh to calculate the cell volumes over the original shape of orography (hereafter the ``true cell volumes").
However we want to avoid the expense of conservative mapping every time step, particularly as we are only interested in the true cell volumes and not in the shape of orography itself.
Therefore, instead of tracking the shape of orography within each cell, we track the true cell volumes by solving the following advection equation for a cell volume:
\begin{eqnarray}
\frac{A^{n+1}V^{n+1}-A^{n}V^{n}}{\Delta t} &=& \sum\limits_{\rm faces}\tilde{A^{n}_{f}}\phi_{m}, \label{aeq}
\end{eqnarray}
where $A$ is the volume adjustment parameter and $V$ is the volume of cells as defined only by their vertices so that $AV$ corresponds to the true cell volumes (Figure \ref{fig:AV}).
To avoid having negative cell volumes, we use a \textit{downwind} value of $A^{n}_{f}$, denoted by $\tilde{A^{n}_f}$, in the right-hand side.
This guarantees that $A$ is always positive as long as the initial value of $A$ is positive at all cells (see appendix A for a proof).
Then we use $AV$ in the advection equation \eqref{discretised_form} instead of $V$ and solve it using a two stage Runge-Kutta method with an off-centring parameter $\alpha$:
\begin{eqnarray}
\frac{A^{n+1}V^{n+1}\rho^{*}-A^{n}V^{n}\rho^{n}}{\Delta t} + (1-\alpha)\sum\limits_{\rm faces}\rho^{n}_{f}(\phi^{n} - \tilde{A^{n}_{f}}\phi_{m})\nonumber \\
+ \alpha\sum\limits_{\rm faces}\rho^{n}_{f}(\phi^{n+1} - \tilde{A^{n}_{f}}\phi_{m}) &=& 0, \label{runge_kutta_1st_stage} \\
\frac{A^{n+1}V^{n+1}\rho^{n+1}-A^{n}V^{n}\rho^{n}}{\Delta t} + (1-\alpha)\sum\limits_{\rm faces}\rho^{n}_{f}(\phi^{n} - \tilde{A^{n}_{f}}\phi_{m})\nonumber \\
+ \alpha\sum\limits_{\rm faces}\rho^{*}_{f}(\phi^{n+1} - \tilde{A^{n}_{f}}\phi_{m}) &=& 0, \label{runge_kutta_2nd_stage}
\end{eqnarray}
where the mesh flux $\phi_m$ is evaluated at time step $n+1/2$ as in Figure \ref{fig:mesh_flux}. In this way, we achieve conservation of both the total of the true cell volume (i.e., the total $AV$) and the total mass relative to the true domain size (i.e., the total $\rho AV$), thereby maintaining uniform fields without the need to track the shape of orography within each cell.

To prove that our scheme preserves a uniform field on a moving mesh, we assume a divergence free velocity field so that for each cell
\begin{eqnarray}
\sum\limits_{\rm faces}\phi = 0 \label{divergence_free}
\end{eqnarray}
at all time steps, and check if the solution stays uniform when the initial condition is uniform.
Given $\rho^{n} \equiv 1$, the equation \eqref{runge_kutta_1st_stage} becomes
\begin{eqnarray}
\frac{A^{n+1}V^{n+1}\rho^{*}-A^{n}V^{n}}{\Delta t} + (1-\alpha)\sum\limits_{\rm faces}(\phi^{n} - \tilde{A^{n}_{f}}\phi_{m})
+ \alpha\sum\limits_{\rm faces}(\phi^{n+1} - \tilde{A^{n}_{f}}\phi_{m}) = 0. \label{runge_kutta_1st_stage_mod}
\end{eqnarray}
Substituting the equations \eqref{aeq} and \eqref{divergence_free} into \eqref{runge_kutta_1st_stage_mod}, we have
\begin{eqnarray}
A^{n+1}V^{n+1}(\rho^{*}-1) = 0.
\end{eqnarray}
As $A^{n+1}V^{n+1} \neq 0$, we obtain $\rho^{*} \equiv 1$. Then the equation \eqref{runge_kutta_2nd_stage} becomes
\begin{eqnarray}
\frac{A^{n+1}V^{n+1}\rho^{n+1}-A^{n}V^{n}}{\Delta t} + (1-\alpha)\sum\limits_{\rm faces}(\phi^{n} - \tilde{A^{n}_{f}}\phi_{m}) + \alpha\sum\limits_{\rm faces}(\phi^{n+1} - \tilde{A^{n}_{f}}\phi_{m}) = 0. \label{runge_kutta_2nd_stage_mod}
\end{eqnarray}
In the same way as above, we obtain $\rho^{n+1} \equiv 1$. Therefore it is proved that the solution stays uniform in a divergence-free velocity field when the initial condition is uniform. In section \ref{results}, we will confirm this result numerically, whereas the model without the volume adjustment suffers from artificial compression and expansion of the fluid in association with the mesh movement over orography.

\subsection{Advection scheme}
\label{advectionScheme}

Section \ref{discretisation} describes the interaction of the discretisation with the moving mesh, and section \ref{colin_parameter} includes the description of a two-stage, second-order Runge-Kutta time stepping scheme as in the equations \eqref{runge_kutta_1st_stage} and \eqref{runge_kutta_2nd_stage}. To complete the discretisation we must specify how face values, $\rho_f$, are calculated from cell values, $\rho$. We use a simple, second-order linear upwind advection scheme without monotonicity constraints. The use of an unbounded advection scheme makes it easier to ensure that the mesh motion over orography does not generate spurious oscillations. The face values are approximated as:
\begin{equation}
\rho_f = \rho_u + \bm{\delta} \cdot \nabla_u \rho
\end{equation}
where $\rho_u$ is the value of $\rho$ in the cell upwind of the face, $\bm{\delta}$ is the vector that goes from the upwind cell centre to the face centre, and $\nabla_u \rho$ is the gradient of $\rho$ calculated in the upwind cell using Gauss' divergence theorem:
\begin{equation}
\nabla_u \rho = \frac{1}{V}\sum\limits_{\text{faces of }u} \tilde{\rho}_f \mathbf{S}_f
\end{equation}
where $\tilde{\rho}_f$ is the values of $\rho$ linearly interpolated from cell centres onto faces and $\mathbf{S}_f$ is the outward pointing vector normal to each face with magnitude equal to the face area (the face area vector).

\section{Results}\label{results}

\subsection{Advection over smooth orography}\label{doubleCosine}
In this section, we present the results of an advection test on a three-dimensional mesh with one layer in the vertical direction.
We use a computational domain with a size of $[-L, L] \times [-L, L] \times [0, H]$, where the domain half-length $L$ and height $H$ are set to 5 km and 1 km, respectively.
The number of cells is $N$ both in the $x$ and $y$ directions.
All boundaries of the domain are considered as a rigid wall.

\begin{figure}[t]
  \centering
  \begin{subfigure}{0.5\hsize}
    \centering
    {\small (a})\\
    \includegraphics[height=50mm, clip, trim=0 40 0 0]{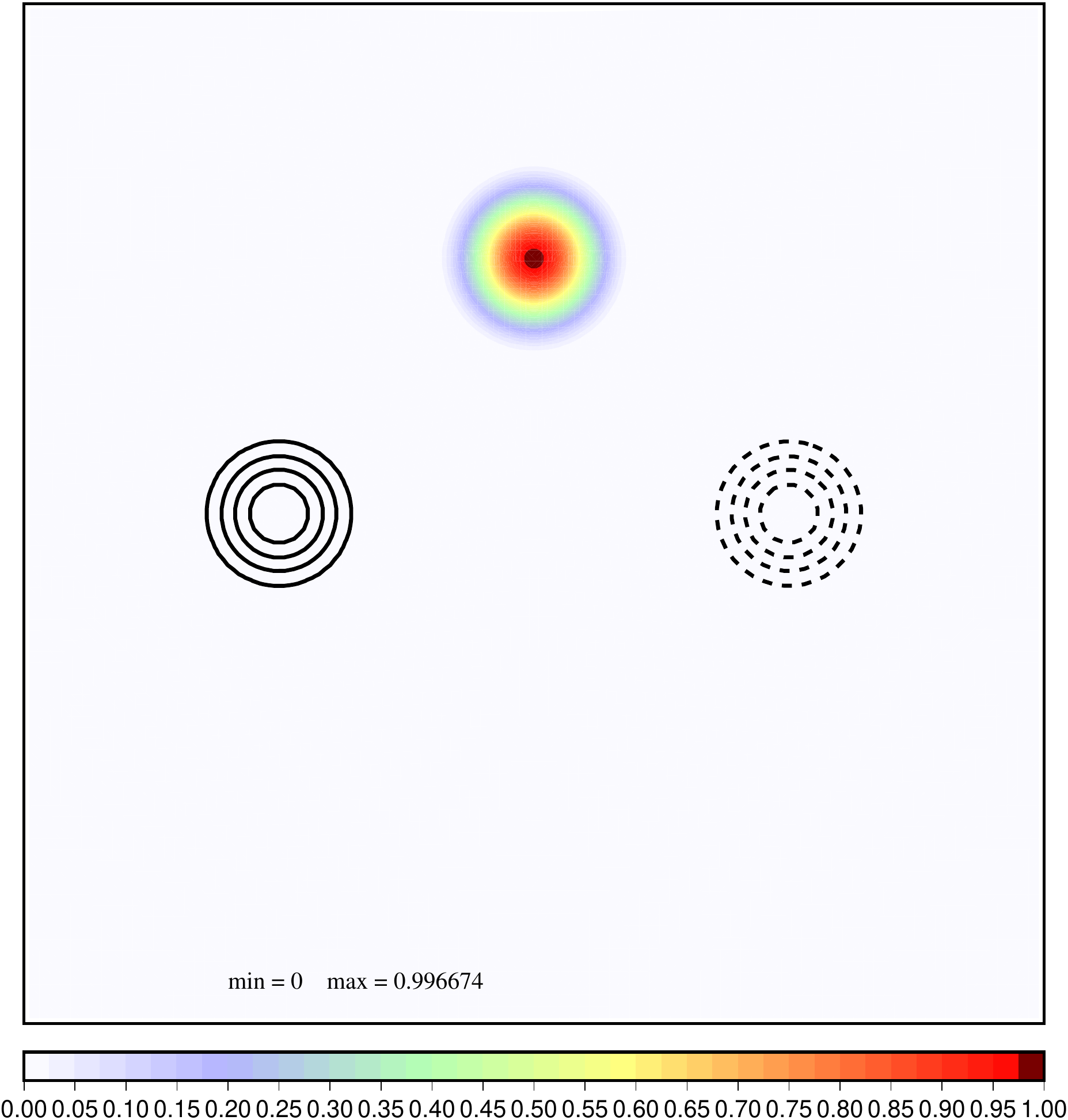} 
    \includegraphics[width=70mm, clip, trim=0 0 0 520]{figures/initial_tracer.pdf} 
  \end{subfigure}
  \begin{subfigure}{0.28\hsize}
    \centering
    {\small (b})
    \includegraphics[height=56mm, clip, trim=0 -70 0 -3]{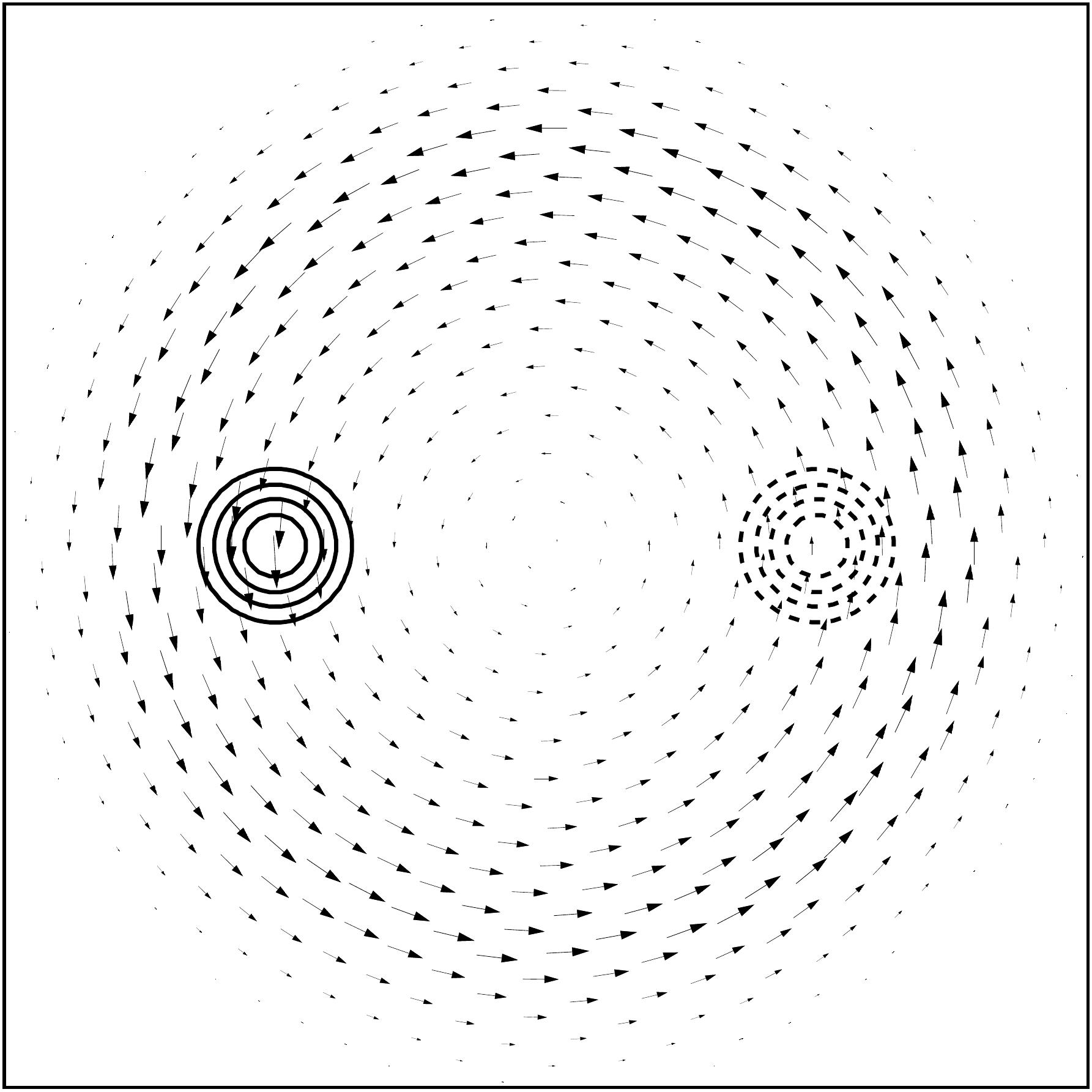} 
  \end{subfigure}
  \caption{Initial conditions of (a) the tracer density and (b) the velocity along with the profile of orography implemented at the bottom of the domain. Colour contours show the amplitudes of the initial tracer density $\rho_{0}$. The arrows represent the initial velocity vector $\MM{u}$. Solid and dashed lines indicate the positive and negative height of orography, respectively, with the contour interval of 100 m.}
  \label{fig:initial_settings}
\end{figure}
\begin{figure}[t]
  \centering
  \begin{subfigure}{0.3\hsize}
    \centering
    {\small (a)}
    \includegraphics[height=48mm]{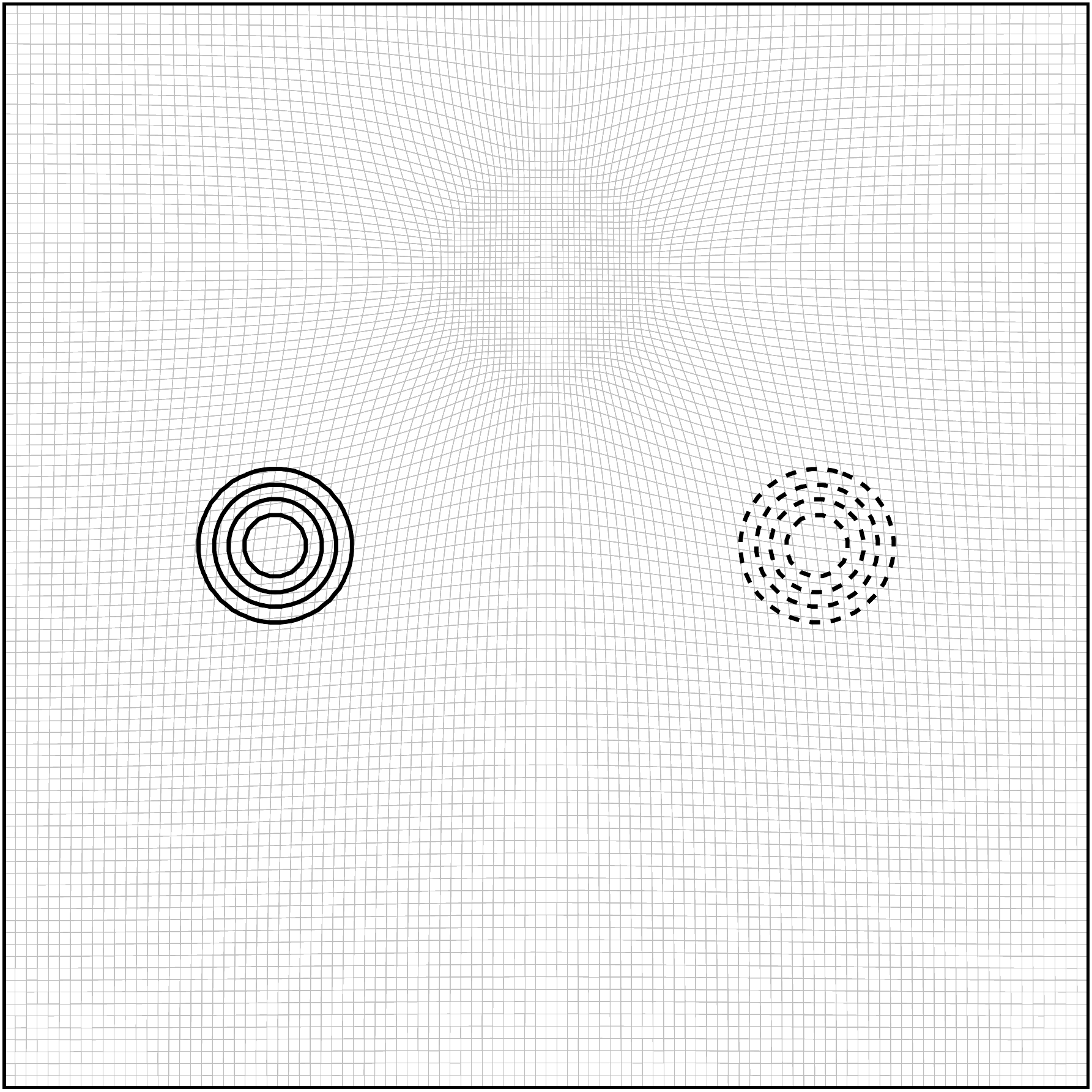} 
  \end{subfigure}
  \begin{subfigure}{0.4\hsize}
    \centering
    {\small (b)}
    \includegraphics[height=40mm]{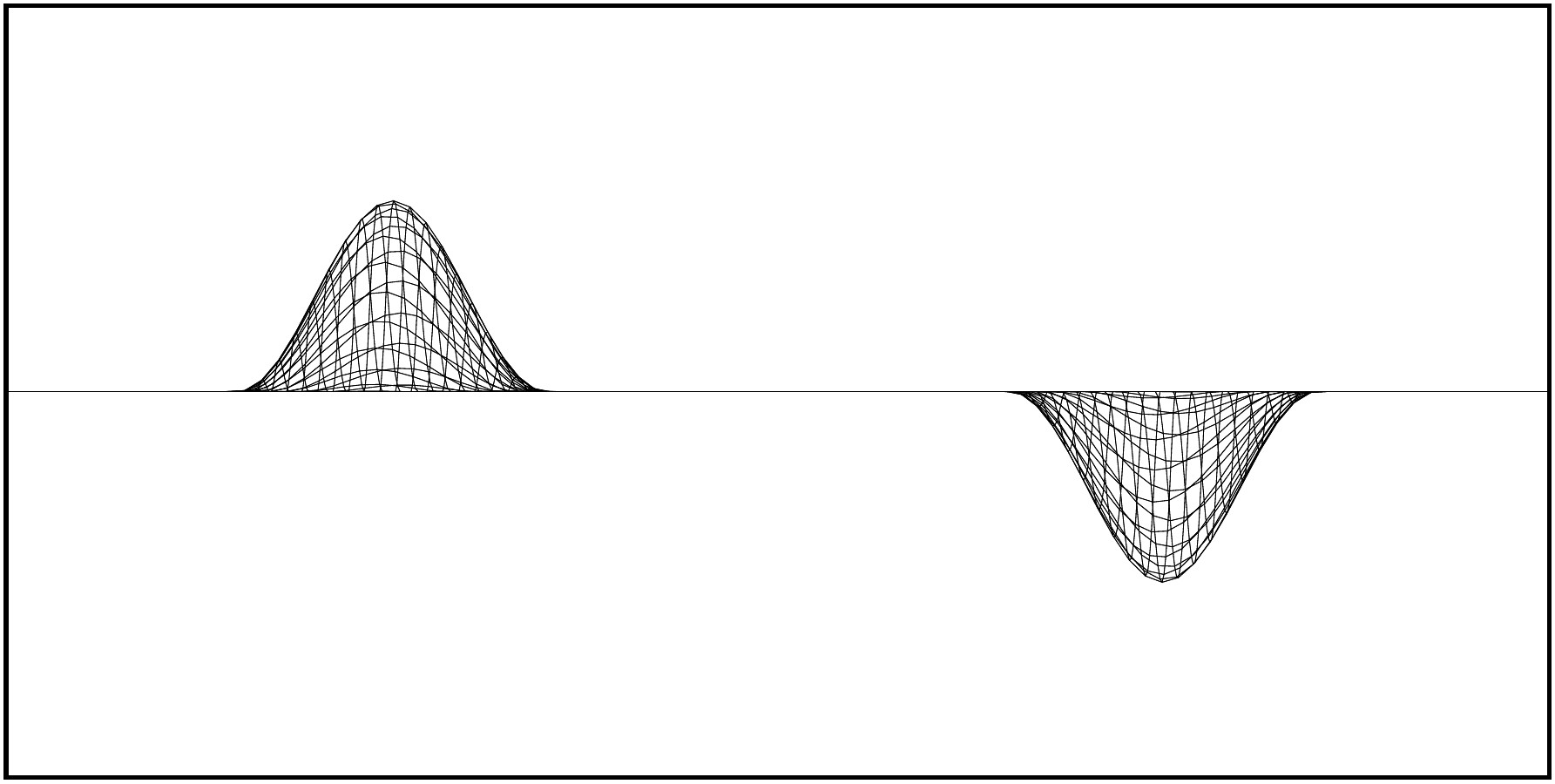} 
  \end{subfigure}
  \caption{The mesh at the initial state. (a) The horizontal x-y slice of the initial mesh at the ground level. The contours show the profile of orography as in Figure \ref{fig:initial_settings}. (b) The vertical x-z slice of the mesh at the ground level through the centre of the hill and valley.}
  \label{fig:mountains}
\end{figure}

\begin{figure}[t]
  \centering
  \begin{subfigure}{0.24\hsize}
    \centering
    {\small (a)}
    \includegraphics[height=42mm, clip, trim=0 40 0 0]{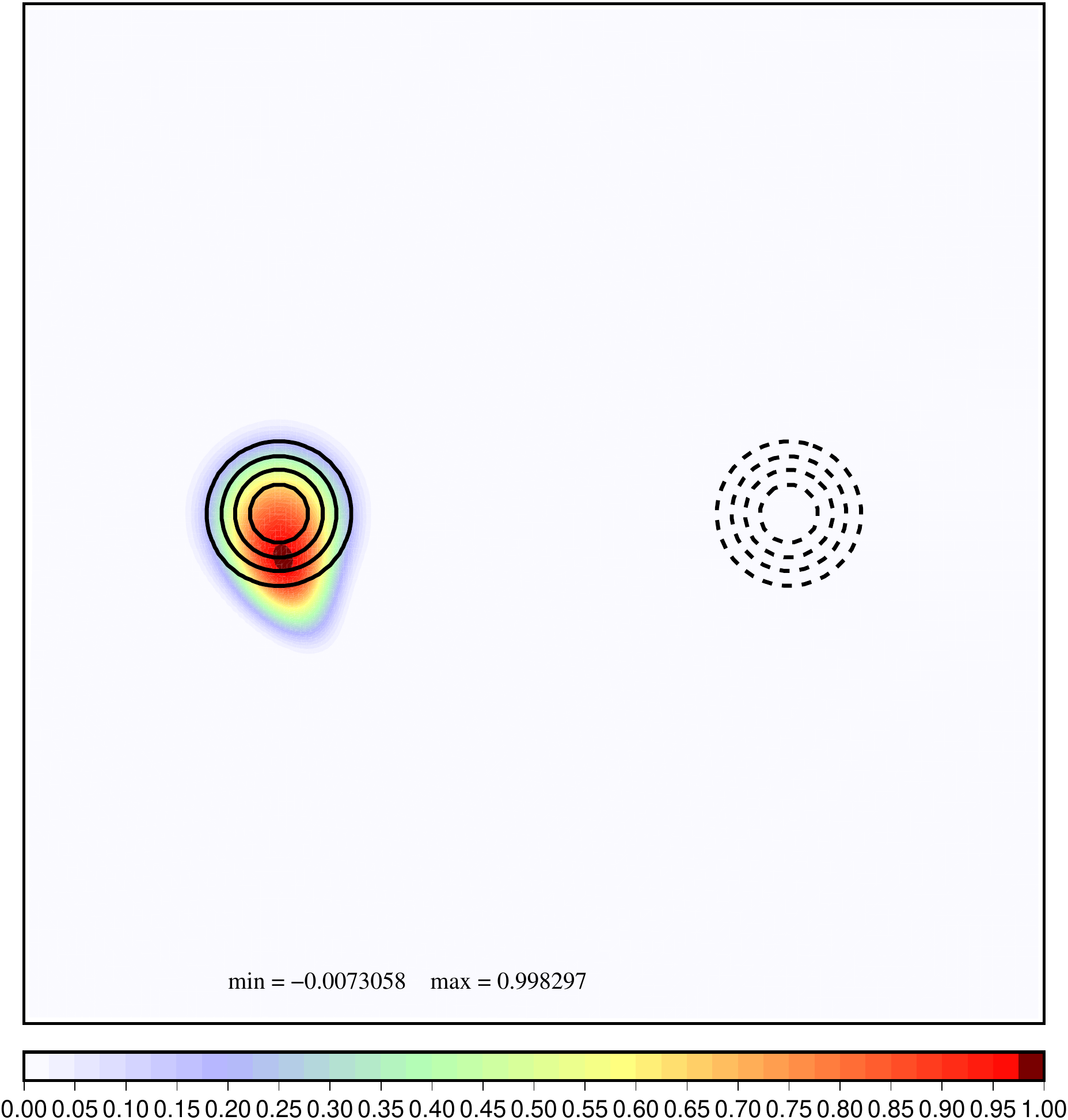} 
  \end{subfigure}
  \begin{subfigure}{0.24\hsize}
    \centering
    {\small (b)}
    \includegraphics[height=42mm, clip, trim=0 40 0 0]{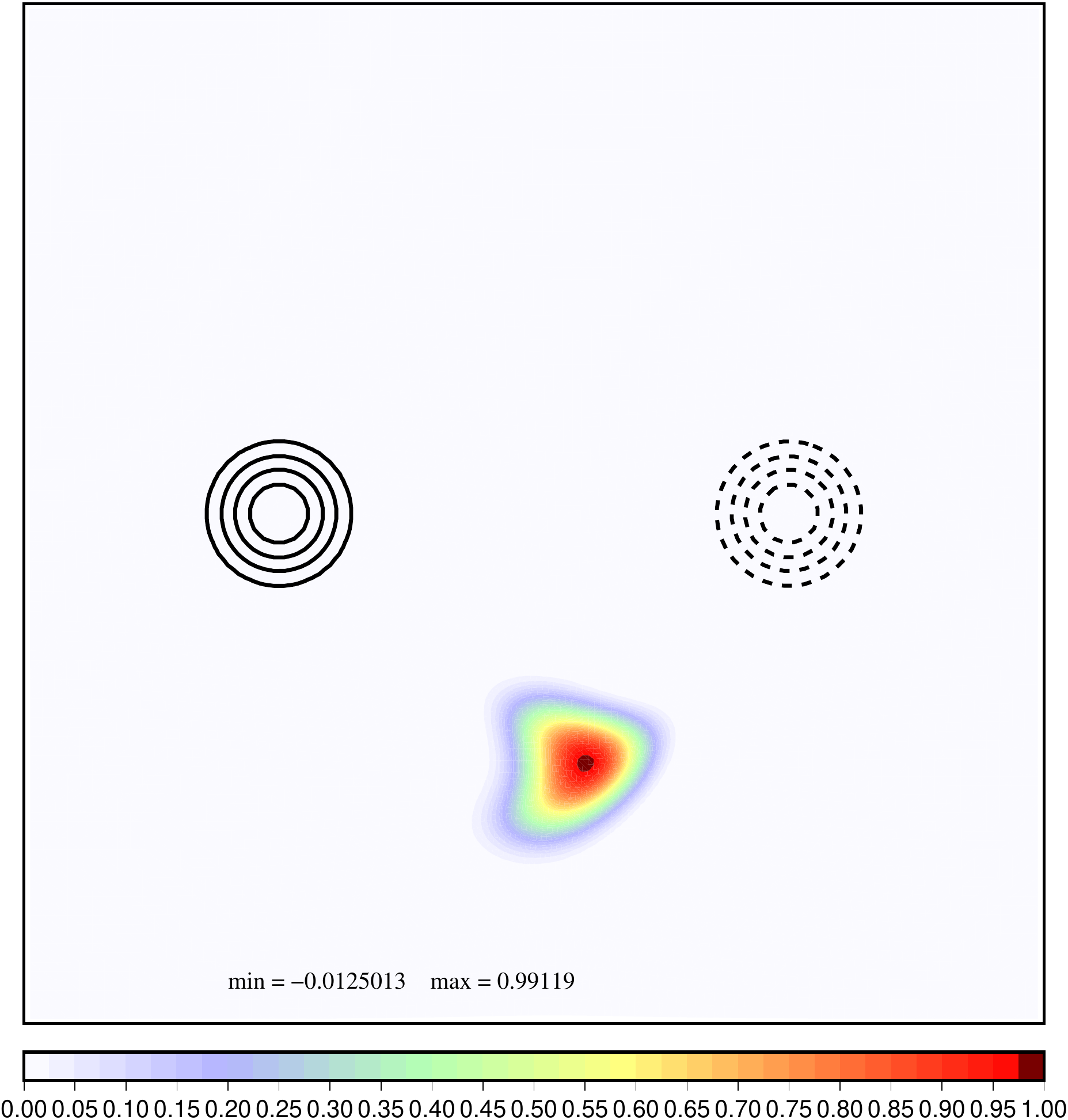} 
  \end{subfigure}
  \begin{subfigure}{0.24\hsize}
    \centering
    {\small (c)}
    \includegraphics[height=42mm, clip, trim=0 40 0 0]{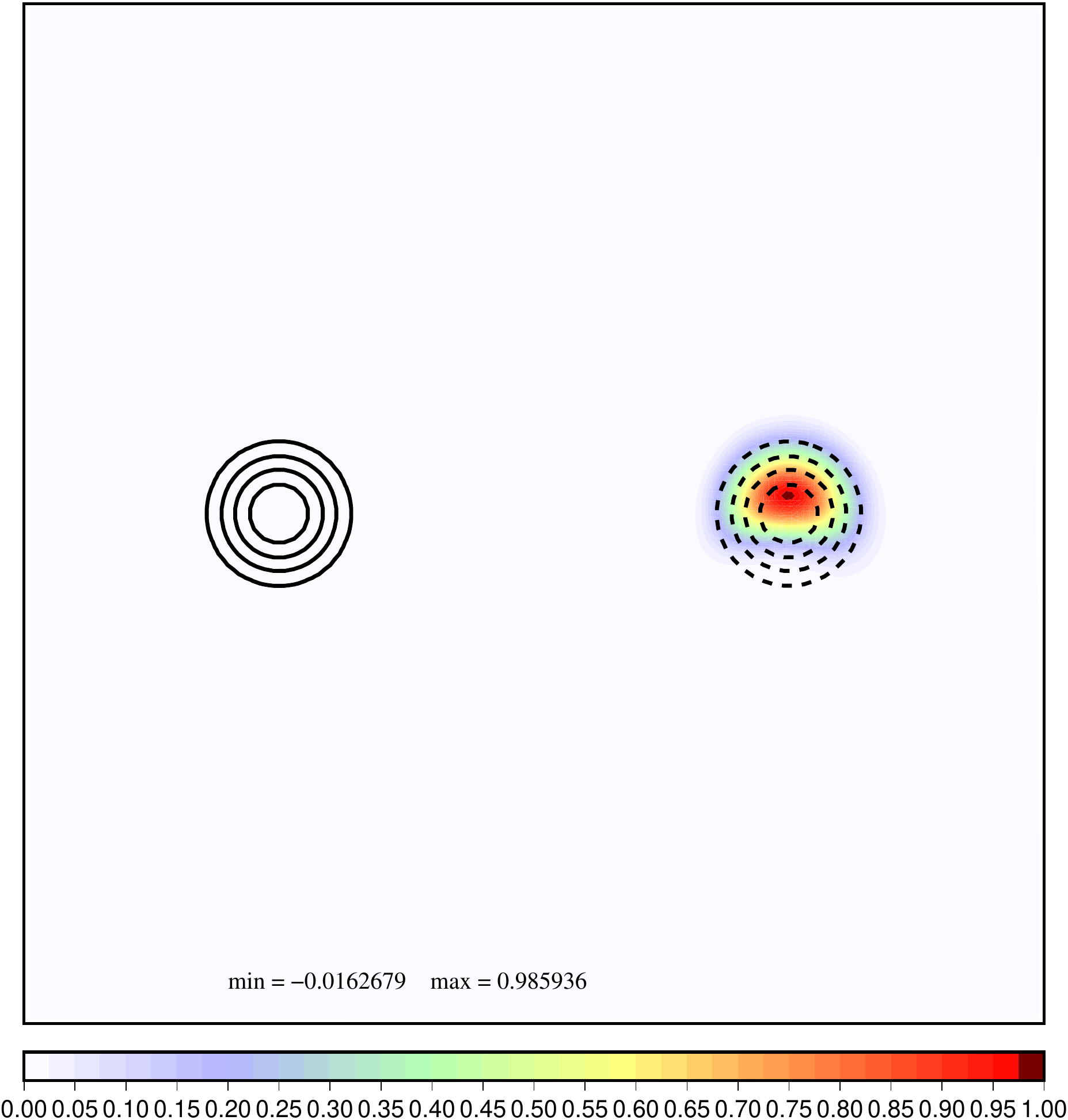} 
  \end{subfigure}
  \begin{subfigure}{0.24\hsize}
    \centering
    {\small (d)}
    \includegraphics[height=42mm, clip, trim=0 40 0 0]{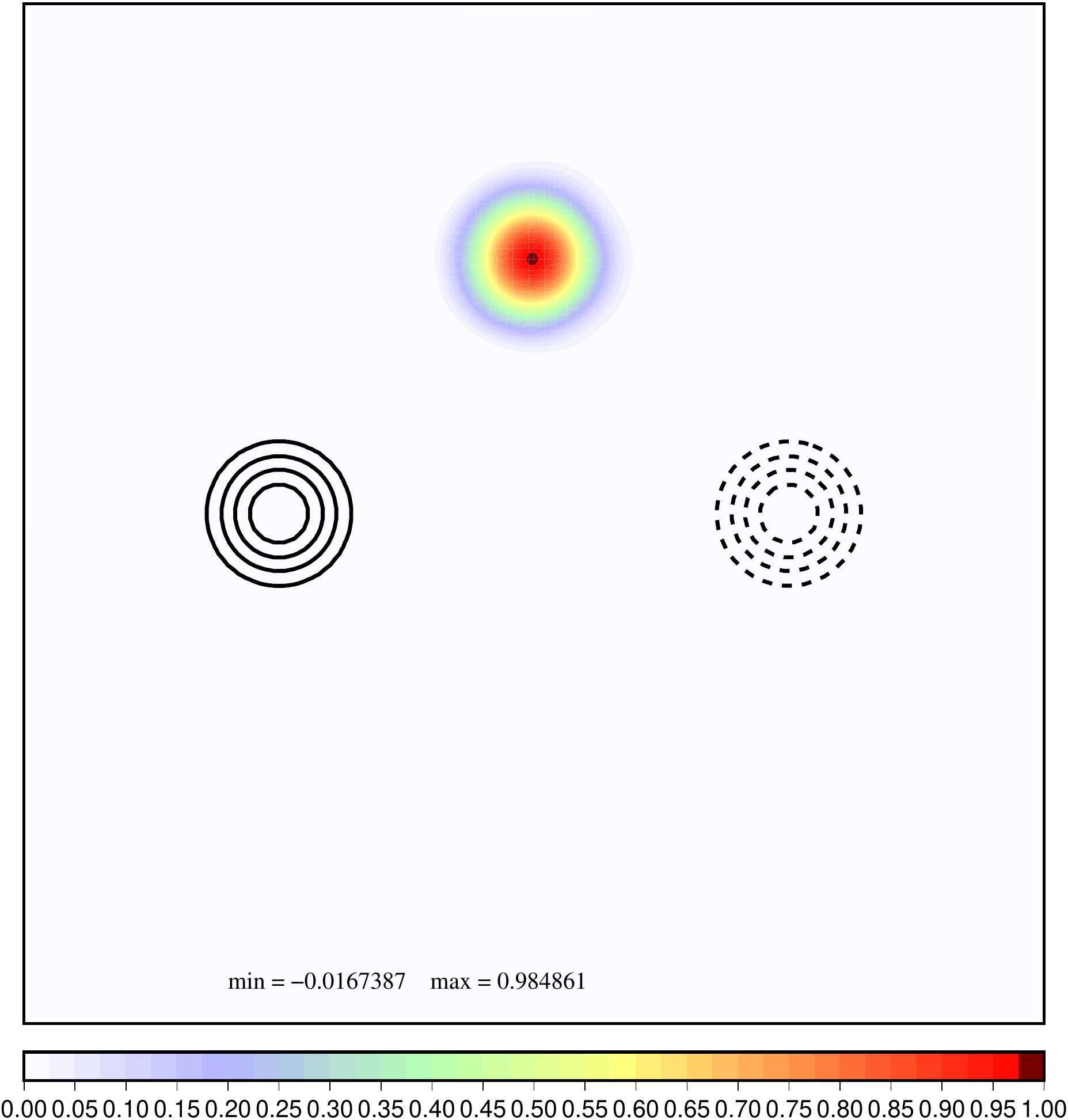} 
  \end{subfigure}
  \begin{subfigure}{0.5\hsize}
    \centering
    \includegraphics[width=70mm, clip, trim=0 0 0 520]{figures/tracer_150.pdf} 
  \end{subfigure}
  \caption{Results of the advection test over smooth orography using a cosine-shaped tracer. Snapshots are taken at $t$ = (a) 150 s, (b) 300 s, (c) 450 s and (d) 600 s. Colour contours show the amplitudes of the tracer density $\rho$. Solid and dashed lines indicate the positive and negative height of orography, respectively, where the contour interval is 100 m.}
  \label{fig:tracer}
\end{figure}
\begin{figure}[h]
  \centering
  \begin{subfigure}{0.24\hsize}
    \centering
    {\small (a)}
    \includegraphics[height=41.5mm]{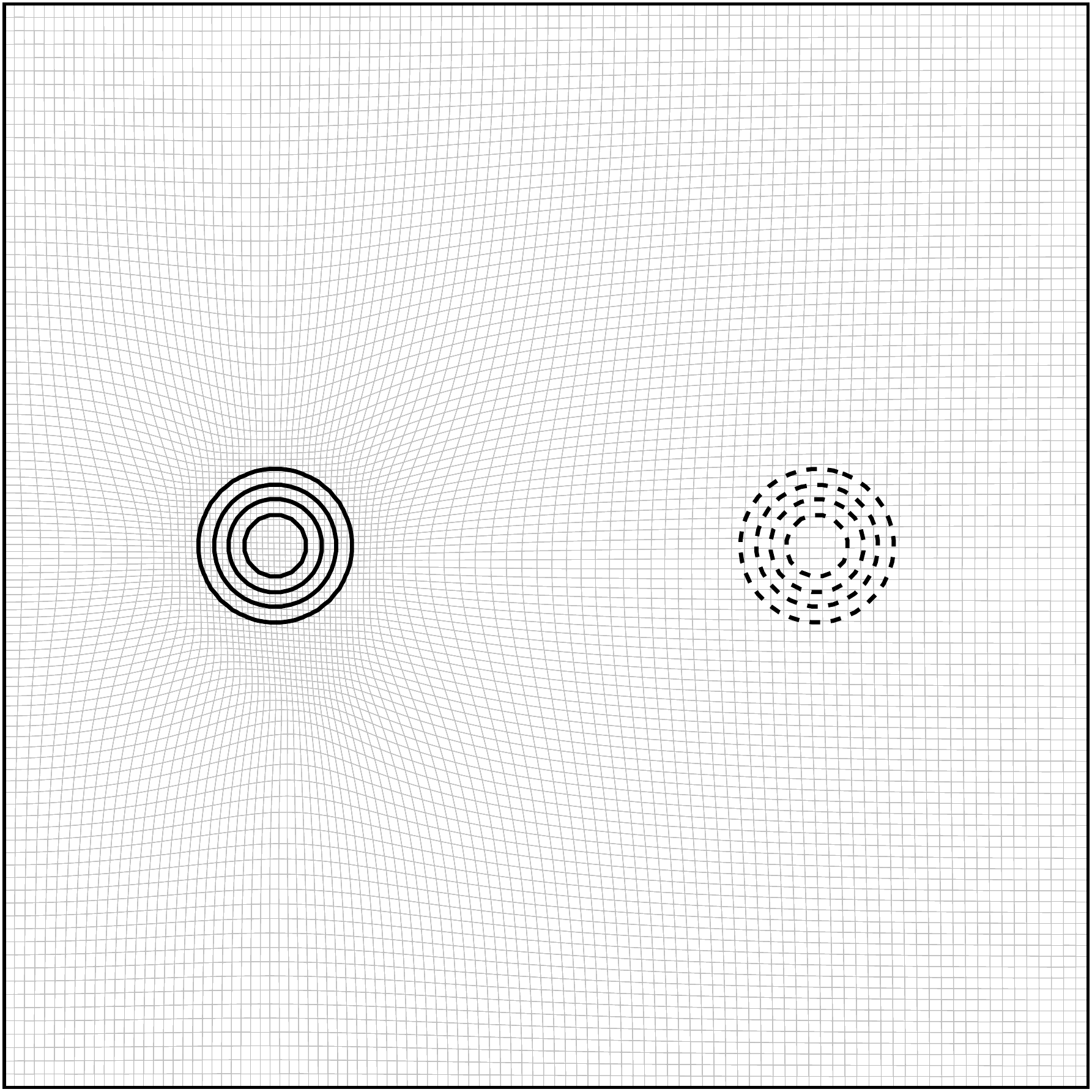} 
  \end{subfigure}
  \begin{subfigure}{0.24\hsize}
    \centering
    {\small (b)}
    \includegraphics[height=41.5mm]{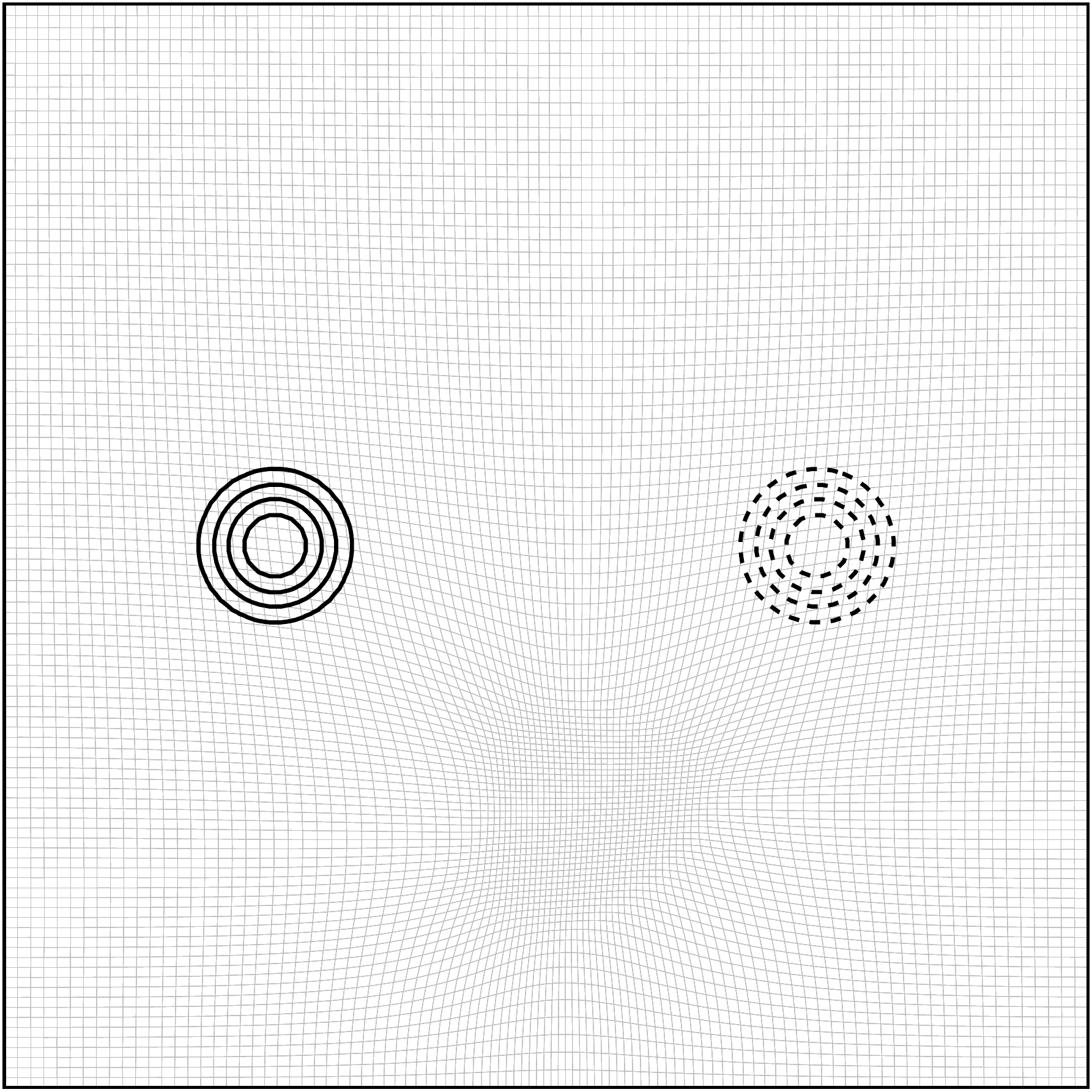} 
  \end{subfigure}
  \begin{subfigure}{0.24\hsize}
    \centering
    {\small (c)}
    \includegraphics[height=41.5mm]{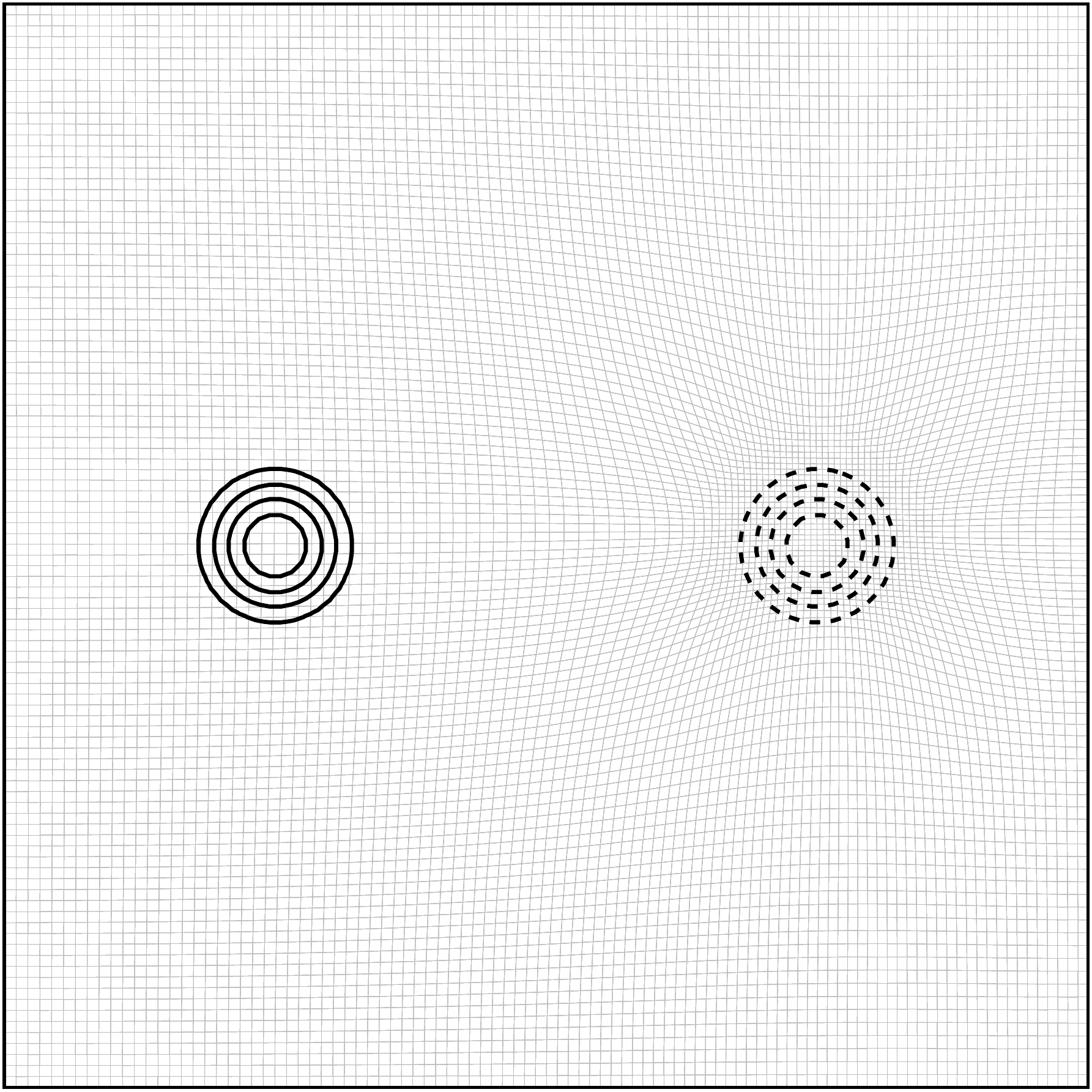} 
  \end{subfigure}
  \begin{subfigure}{0.24\hsize}
    \centering
    {\small (d)}
    \includegraphics[height=41.5mm]{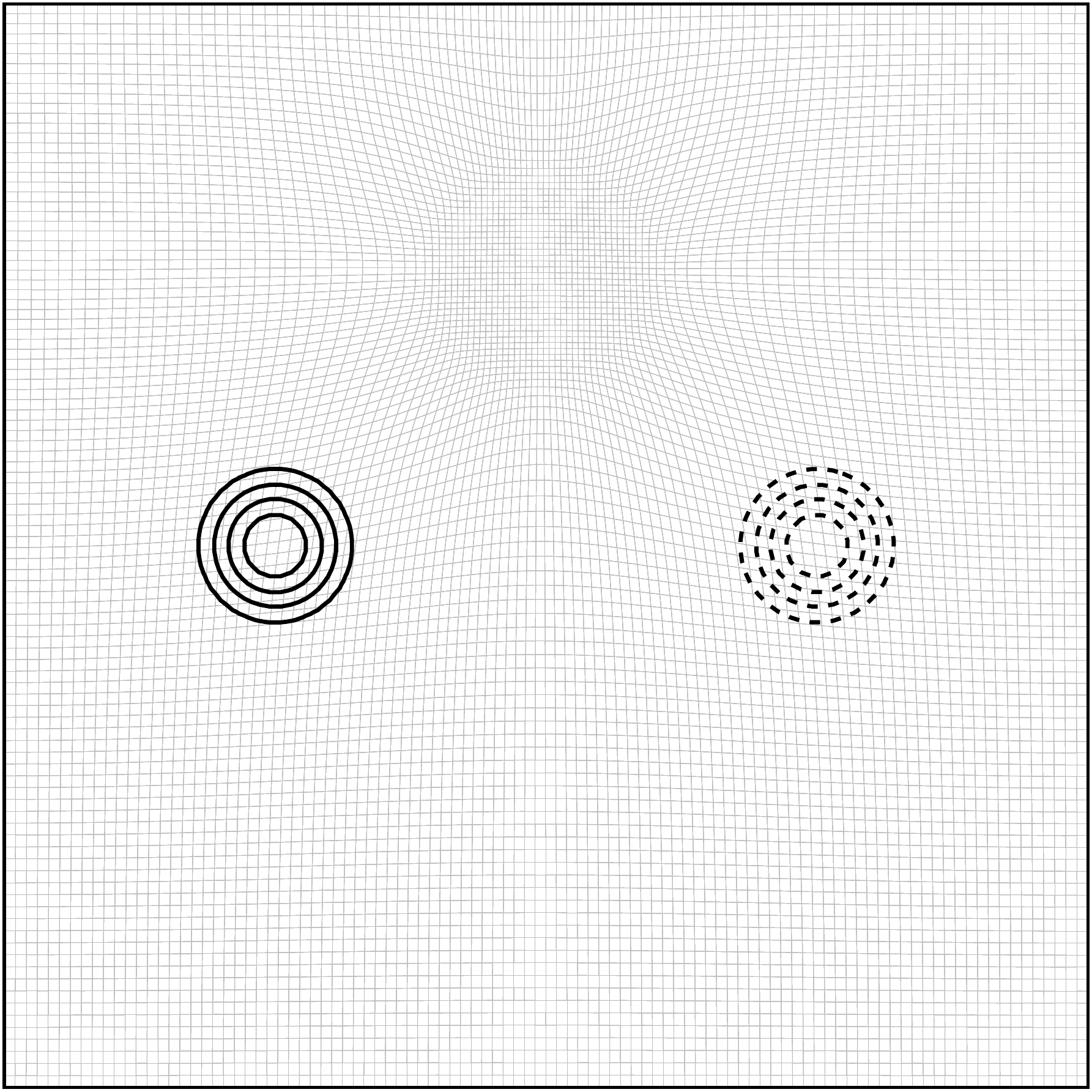} 
  \end{subfigure}
  \caption{Snapshots of the moving mesh at $t$ = (a) 150 s, (b) 300 s, (c) 450 s and (d) 600 s in the advection test over smooth orography using a cosine-shaped tracer. The contours in the background show the profile of orography as in Figure \ref{fig:tracer}.}
  \label{fig:mesh}
\end{figure}
\begin{figure}[h]
  \begin{subfigure}{0.49\hsize}
  \leftline{\small (a)}
  \centering
  \includegraphics[height=60mm]{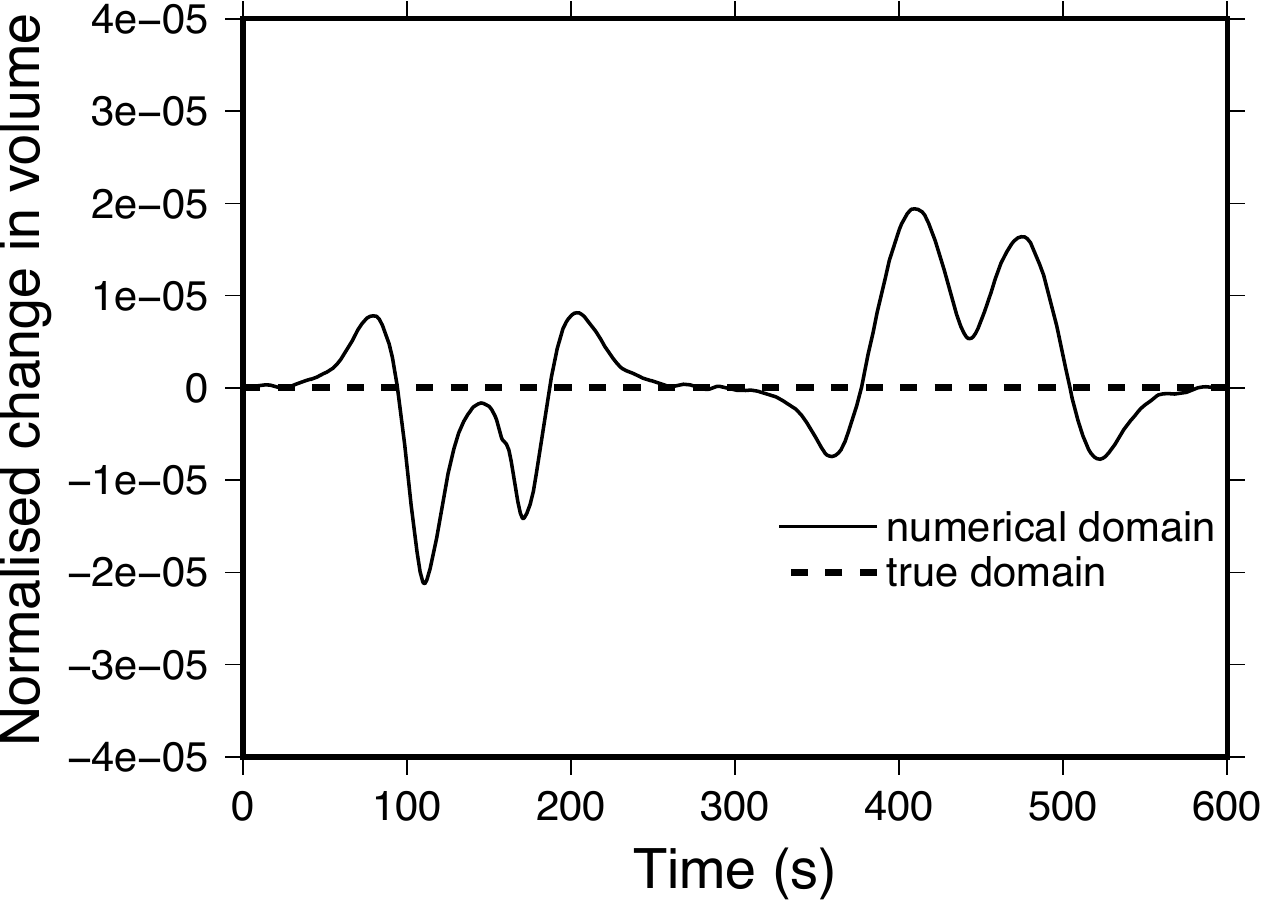} 
  \end{subfigure}
  \begin{subfigure}{0.49\hsize}
  \leftline{\small (b)}
  \centering
  \includegraphics[height=60mm]{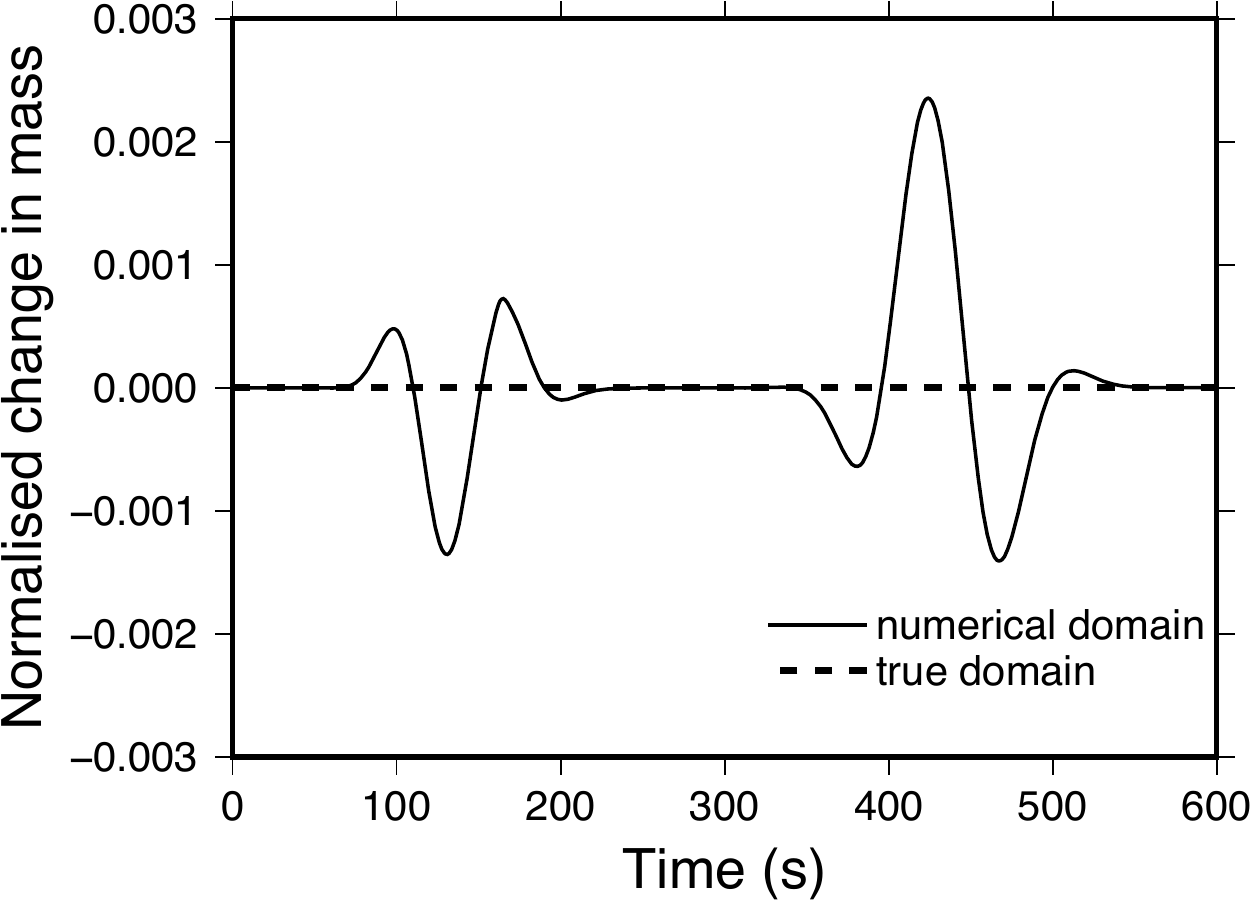} 
  \end{subfigure}
  \caption{Conservation of total cell volumes and total mass in the advection test over smooth orography using a cosine-shaped tracer. (a) Solid line shows the error in the total $V$ relative to the initial domain size, and dashed line shows that of $AV$. (b) Solid line shows the error in the total mass calculated in the numerical domain (i.e. the sum of all $\rho V$) relative to the initial total mass, and dashed line shows the error in the total mass calculated in the true domain (the sum of all $\rho AV$) relative to the initial total mass.}
  \label{fig:diagnostics}
\end{figure}
\begin{figure}[h]
  \centering
  \includegraphics[height=60mm]{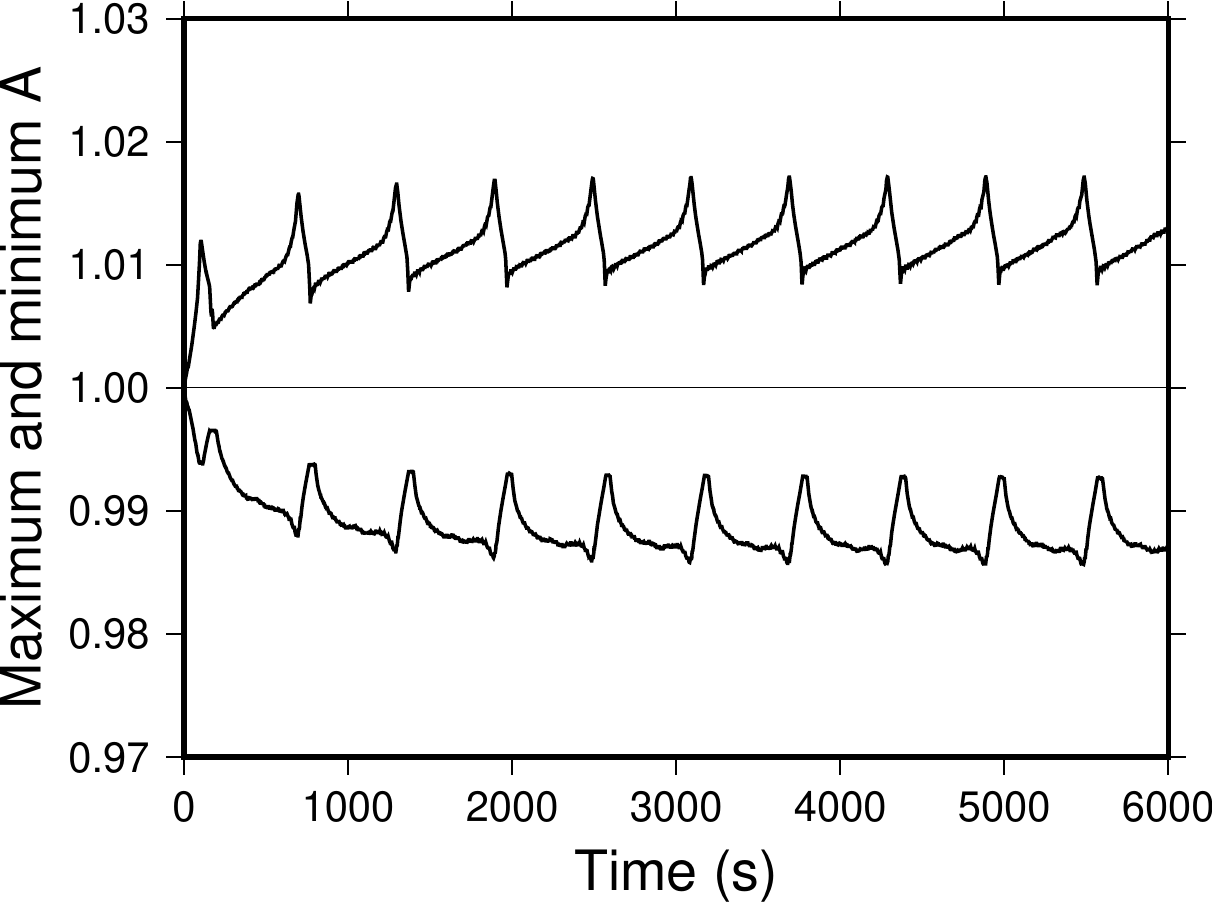} 
  \caption{Time evolution of the volume adjustment parameter $A$ in the advection test over smooth orography using a cosine-shaped tracer. Upper and lower lines show the maximum and minimum values of $A$ in the domain, respectively, for the period of 10 revolutions.}
  \label{fig:Aminmax}
\end{figure}

A bubble of tracer is transported by a solid body rotating velocity field.
Figure \ref{fig:initial_settings} shows diagrams of the tracer density and the velocity of the initial state.
The initial tracer density $\rho_0$ is defined as
\begin{empheq}[left={\rho_0(x, y) =\empheqlbrace}]{alignat=2}
& \frac{1}{2}\left[1+\cos\left(\frac{\pi r_t}{R_t}\right)\right] & \quad (r_t \le R_t) \\
& 0 & \quad (r_t > R_t)
\end{empheq}
where the tracer radius $R_t = L/5$, and the distance to the centre of the tracer,
\begin{eqnarray}
r_t &=& |\mathbf{x} - \mathbf{x}_t| = 
\sqrt{(x-x_t)^{2}+(y-y_t)^{2}},
\end{eqnarray}
with the centre of the tracer initially at $\mathbf{x}_t = (0, L/2)$.
The non-divergent velocity field can be written in terms of a stream function $\psi$ as
\begin{eqnarray}
u = - \frac{\partial \psi}{\partial y}, \; \; v = \frac{\partial \psi}{\partial x},
\end{eqnarray}
and $w = 0$. We use $\psi$ that yields a velocity field which rotates around the centre of the domain and decays linearly to zero before it reaches the boundaries:
\begin{empheq}[left={\psi(x,y) =\empheqlbrace}]{alignat=3}
&\Omega\  r_v^2 & \quad (r_v \le R_i) \\
&\Omega\ R_i \left\{ R_i + (r_v-R_i)\left(\frac{R_o - r_v}{R_o-R_i} + 1\right)\right\}  & \quad (R_i < r_v \le R_o) \\
&\Omega\ R_i R_o & (r_v > R_o)
\end{empheq}
where $r_v$ is the distance to the centre of the domain. The inner radius is set to $R_i=0.76 L$ so that the tracer is separate from the sheared velocity and the outer radius is set to $R_o=L$ so that the velocity is zero at the boundary. The angular velocity is given by $\Omega=\pi/600$ s$^{-1}$ so that the tracer is transported counterclockwise and reaches its initial position after 600 seconds.
Note that the velocity field is recalculated after every time step so that it is always non-divergent on a moving mesh.

The tracer passes over a hill and a valley, as shown in Figure \ref{fig:mountains}, with surface height given by
\begin{empheq}[left={h(x, y) =\empheqlbrace}]{alignat=2}
& \frac{h_\mathrm{max}}{2}\left[1+\cos\left(\frac{\pi r_h}{a}\right)\right] & \quad (r_h \le a) \\
& \frac{h_\mathrm{min}}{2}\left[1+\cos\left(\frac{\pi r_v}{a}\right)\right] & \quad (r_v \le a) \\
& 0 & \quad (r_h > a \,\, \mathrm{and} \,\, r_v > a)
\end{empheq}
where $a = L/5$ is the orography radius, $h_\mathrm{max} = 500$ m is the height at the centre of the hill and that of the valley $h_\mathrm{min} = - 500$ m. The distance to the centre of the hill $r_h = |\mathbf{x} - \mathbf{x}_h|$ with $\mathbf{x}_h=(-L/2, 0)$, and the distance to the centre of the valley $r_v = |\mathbf{x} - \mathbf{x}_v|$ with $\mathbf{x}_v=(L/2, 0)$.

First we present results from the control run, where the number of cells $N =$ 100 and the time step $\Delta t$ = 0.5 s are used.  
Figures \ref{fig:tracer} and \ref{fig:mesh} show four snapshots of the tracer density and the moving mesh, respectively, at t = 150 s, 300 s, 450 s and 600 s.
As the velocity field is non-divergent, the tracer is accelerated over the hill and decelerated over the valley, returning to its original shape and position after 600 s.
The mesh successfully tracks the tracer as it moves and changes its shape, without tangling.
As described in section \ref{moving_mesh}, the monitor function is chosen here so that the cell areas are a factor of 4 smaller in regions where the second derivatives of the tracer density is highest compared with the regions of lowest second derivatives (see appendix \ref{secn:MAsolution:monitor} for details).
\begin{figure}[t]
  \centering
  \begin{subfigure}{0.3\hsize}
    \centering
    {\small (a)}
    \includegraphics[height=52mm, clip, trim=0 40 0 0]{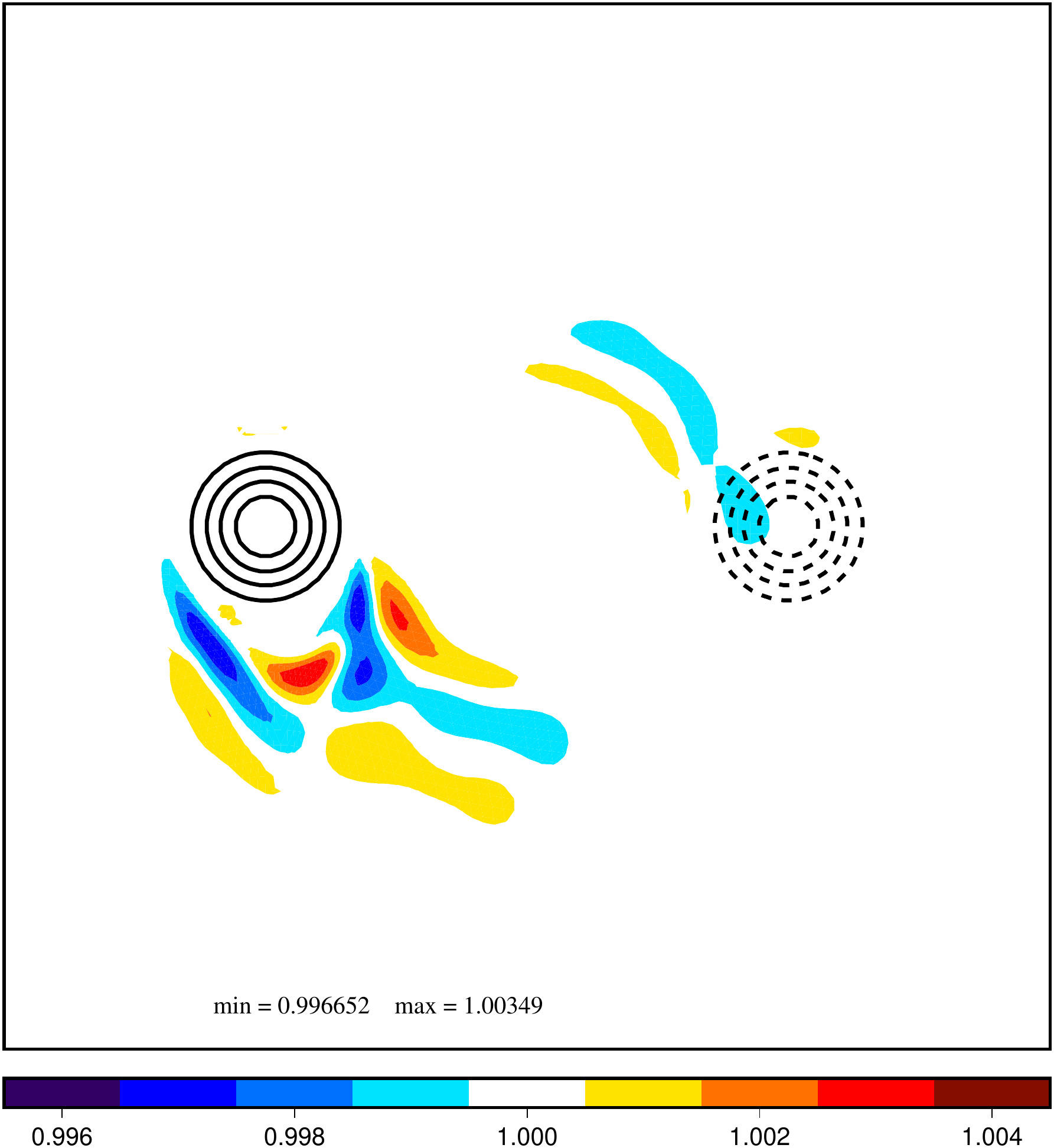} 
  \end{subfigure}
  \begin{subfigure}{0.3\hsize}
    \centering
    {\small (b)}
    \includegraphics[height=52mm, clip, trim=0 40 0 0]{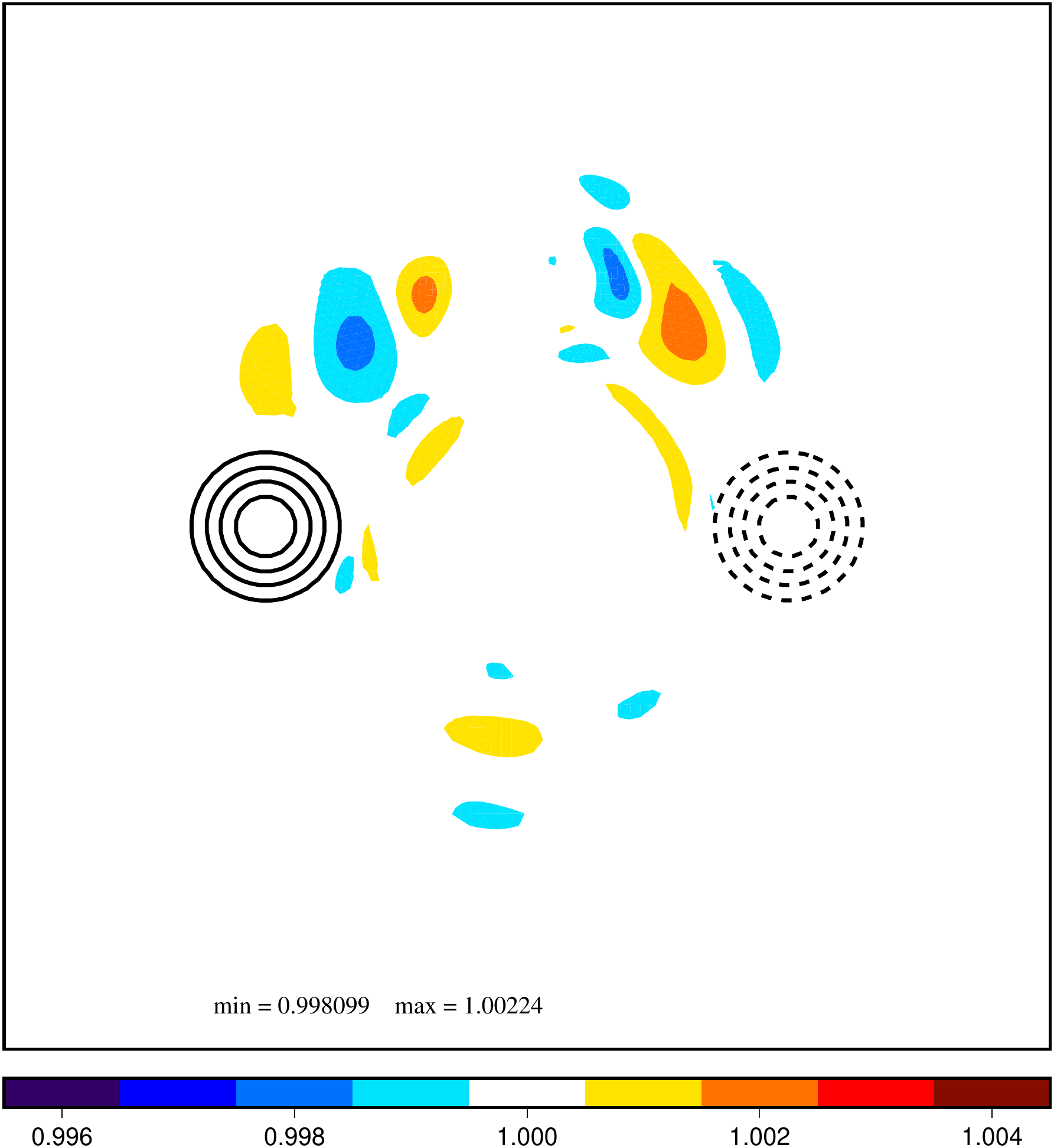} 
  \end{subfigure}
  \begin{subfigure}{0.3\hsize}
    \centering
    {\small (c)}
    \includegraphics[height=52mm, clip, trim=0 40 0 0]{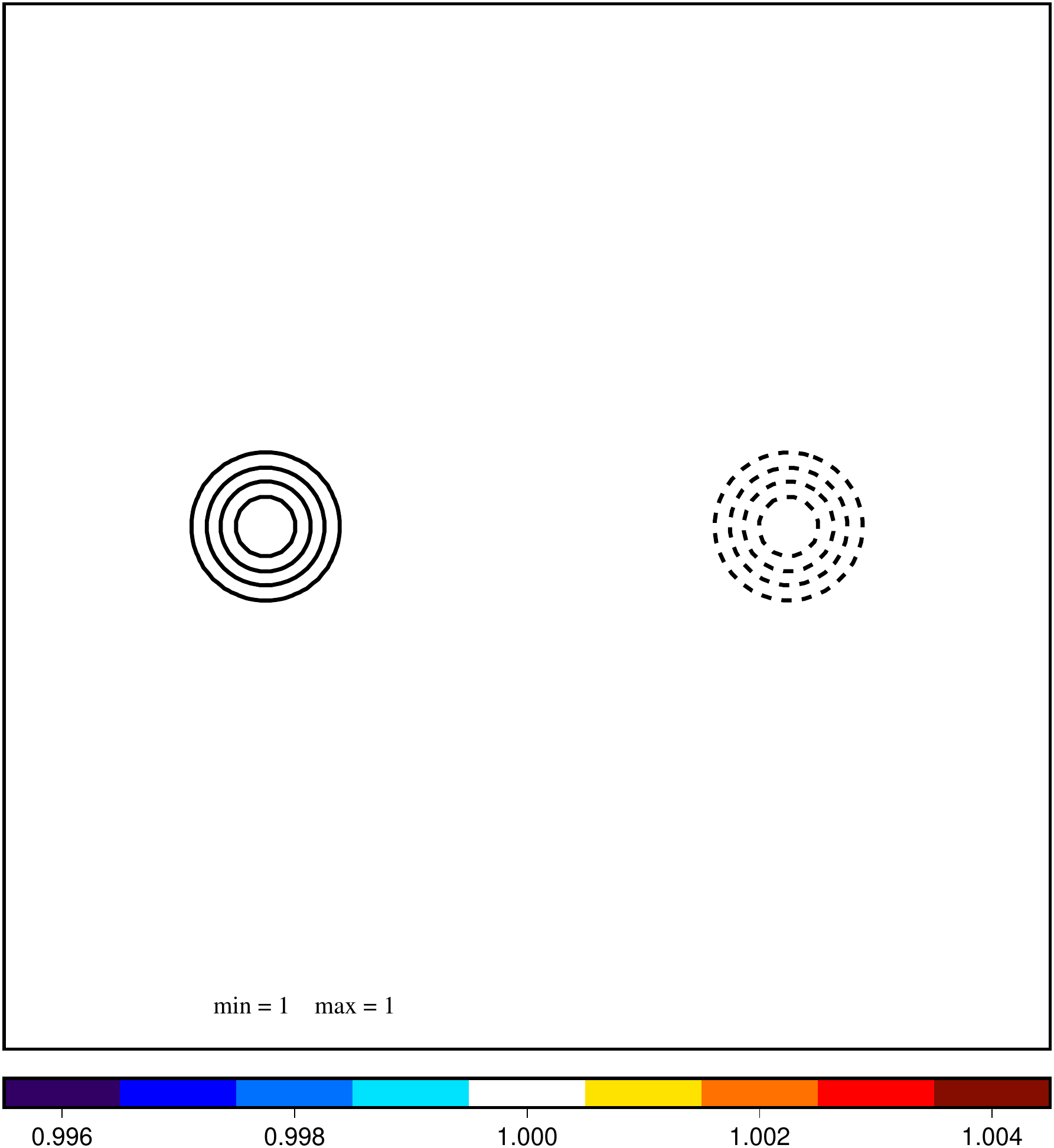} 
  \end{subfigure}
  \centering
  \begin{subfigure}{0.5\hsize}
    \centering
    \includegraphics[width=70mm, clip, trim=0 0 0 520]{figures/uniT_150.pdf} 
  \end{subfigure}
  \caption{Results of the advection test over smooth orography when a uniform density field is used as the initial condition. Figures (a) and (b) show the density fields at t = 150 s and 600 s, respectively, when the volume adjustment parameter $A$ is not used in the model. Figure (c) shows the density field at t = 600 s when $A$ is used to adjust the cell volumes. Colour contours show the amplitudes of the tracer density $\rho$. Solid and dashed lines indicate the positive and negative height of orography, respectively, where the contour interval is 100 m.}
  \label{fig:uniT}
\end{figure}
\begin{figure}[t]
  \centering
  \begin{subfigure}{0.3\hsize}
    \centering
    {\small (a)}
    \includegraphics[height=52mm, clip, trim=0 40 0 0]{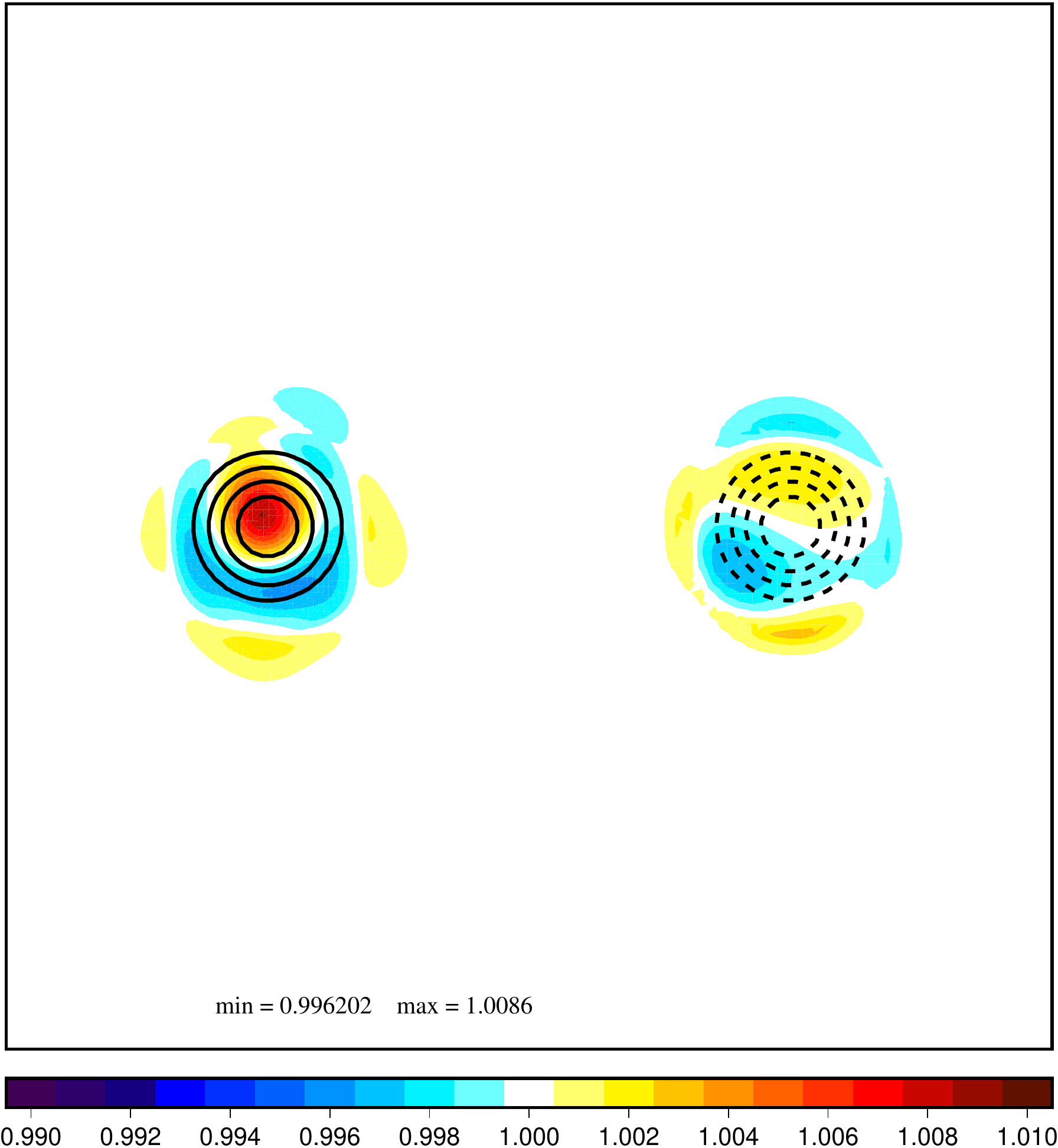} 
  \end{subfigure}
  \begin{subfigure}{0.3\hsize}
    \centering
    {\small (b)}
    \includegraphics[height=52mm, clip, trim=0 40 0 0]{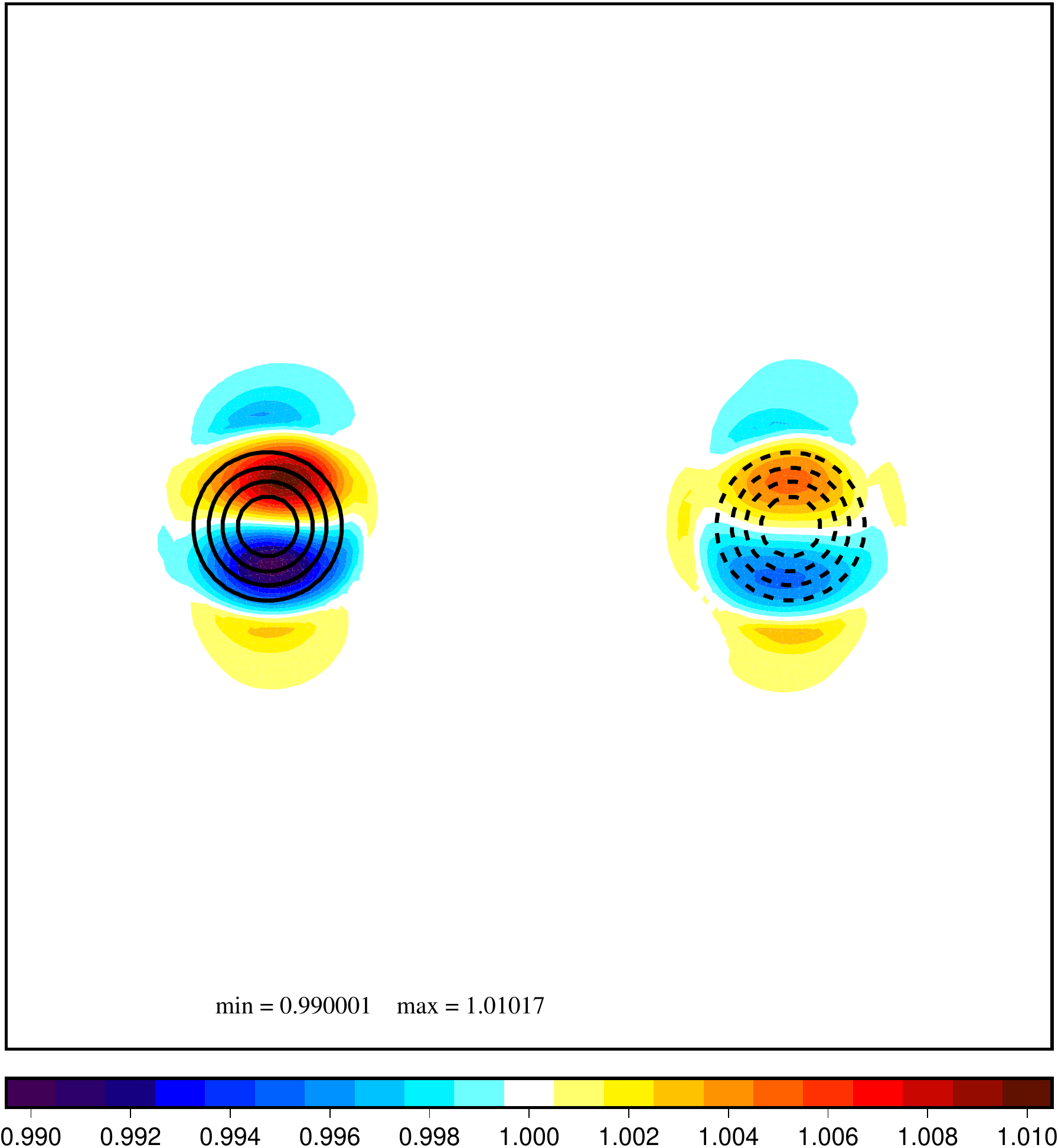} 
  \end{subfigure}
  \begin{subfigure}{0.5\hsize}
    \centering
    \includegraphics[width=70mm, clip, trim=0 0 0 520]{figures/A_150.pdf} 
  \end{subfigure}
  \caption{Snapshots of the $A$ field at $t$ = (a) 150 s and (b) 600 s in the advection test over smooth orography when a uniform density field is used as the initial condition. Colour contours show the amplitudes of the volume adjustment parameter $A$. Solid and dashed contours show the profile of orography as in Figure \ref{fig:uniT}.}
  \label{fig:A}
\end{figure}

Figure \ref{fig:diagnostics} shows the conservation errors in the total cell volumes and the total mass.
While the numerical domain size (the total $V$) changes as the mesh moves, the true domain size (the total $AV$) stays constant (Figure \ref{fig:diagnostics}a).
At the same time, the model conserves the total mass relative to the true domain size (the total $\rho AV$), as shown in Figure \ref{fig:diagnostics}b, achieving exact local conservation and maintenance of uniform fields at the bottom boundary.
Figure \ref{fig:Aminmax} shows the variation of maximum and minimum values of $A$ in the domain for the period of 10 revolutions.
It demonstrates that the first-order downwind differencing of $A$ ensures that $A$ is always positive, thereby the model may not have negative cell volumes. 

To demonstrate that the model maintains uniform fields as mesh moves over orography, we repeat the experiment using a uniform field $\rho_0 \equiv 1$ instead of using the cosine-shaped tracer.
Here we compare the results with and without the use of the volume adjustment parameter $A$.
Figures \ref{fig:uniT}a and \ref{fig:uniT}b show the density fields at $t =$ 150 s and $t =$ 600 s, respectively, when $A$ is not used in the model.
The results show evidence of artificial compression and expansion of the model fluid in association with the mesh movement over orography.
On the other hand, the density field stays constant at 1 throughout the period of 600 s when $A$ is used to adjust the cell volumes (Figure \ref{fig:uniT}c).
Figures \ref{fig:A}a and \ref{fig:A}b show snapshots of the $A$ field at $t =$ 150 s and 600 s, respectively. These results demonstrate that $A$ successfully tracks the changes in the cell volumes over orography and adjusts the cell volumes so that the total of the true cell volumes as well as the total mass in the true domain is conserved, thereby maintaining uniform fields as the mesh moves over orography. 
\begin{figure}[t]
  \centering
  \includegraphics[height=70mm]{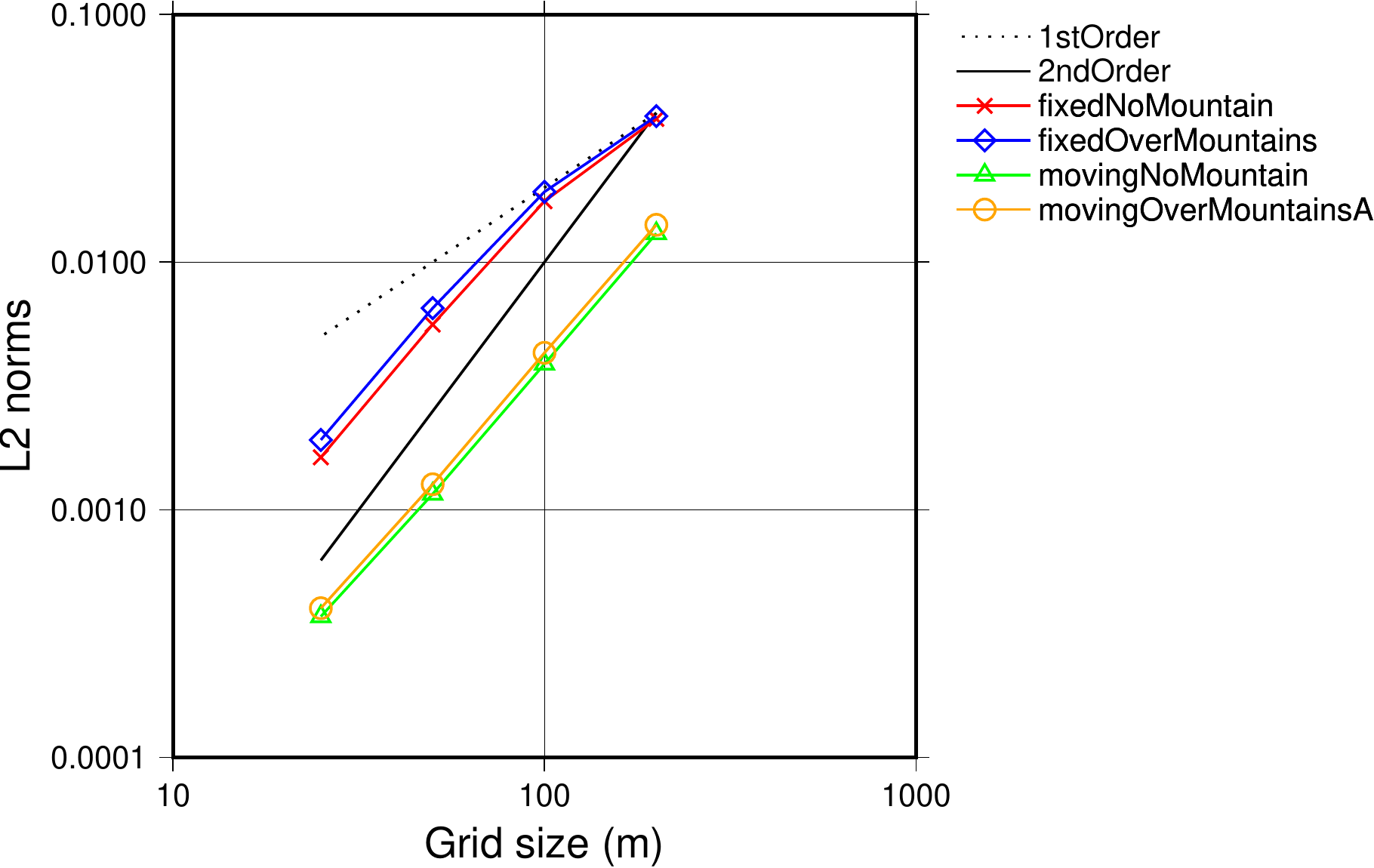} 
  \caption{The $L_2$ norm of the error in $\rho$ versus the grid size for the advection test over smooth orography using a cosine-shaped tracer. Blue and red lines show the error plots calculated on fixed uniform meshes with and without the orography at the bottom boundary, respectively. Orange and green lines show the error plots calculated on moving meshes with and without the orography at the bottom boundary, respectively. Black solid and dashed lines show the theoretical first- and second-order convergence.}
  \label{fig:l2norms}
\end{figure}

Finally, Figure \ref{fig:l2norms} shows a log-log plot of the $L_2$ norm of the errors in $\rho$ versus the grid size after one complete revolution, alongside the theoretical first- and second-order convergence rates.
Time steps of $\Delta t =$ 1 s, 0.25 s and 0.125 s are used for the cases of $N =$ 50, 200 and 400, respectively.
The same cosine-shaped tracer as in the control run is used here at all resolutions.
At each resolution, errors are calculated both on the fixed uniform mesh and the moving mesh, with and without the orography at the bottom surface.
The model shows the convergence rate of 1.70 on the uniform mesh over orography and that of 1.58 on the moving mesh over orography.
It is shown that the implementation of orography doesn't affect the order of convergence both on the uniform mesh and the moving mesh.
The convergence rate on the moving mesh could be improved by optimising the monitor function or using a higher-order advection scheme, but it is outside the scope of this paper.

\subsection{Advection over steep orography}\label{doubleCylinder}

In the previous section, we performed a tracer advection test and showed that our scheme maintains a uniform field on a moving mesh over orography, using smooth hill and valley as shown in Figure \ref{fig:mountains} as a sample orography.
In this section, we demonstrate the importance of the maintenance of uniform fields by repeating the tracer advection test using a pair of cylinder-shaped hill and valley which has steep cliffs on the sides, with surface height given by
\begin{empheq}[left={h(x, y) =\empheqlbrace}]{alignat=2}
& h_c & \quad (r_h \le a) \\
& -h_c & \quad (r_v \le a) \\
& 0 & \quad (r_h > a \,\, \mathrm{and} \,\, r_v > a)
\end{empheq}
where $h_c = 500$ m.
All the other simulation setup is the same as that of the control run in the previous section.
We run the model with and without the use of the volume adjustment parameter $A$ and compare the results.
\begin{figure}[t]
  \centering
  \begin{subfigure}{0.3\hsize}
    \centering
    {\small (a)}
    \includegraphics[height=52mm, clip, trim=0 40 0 0]{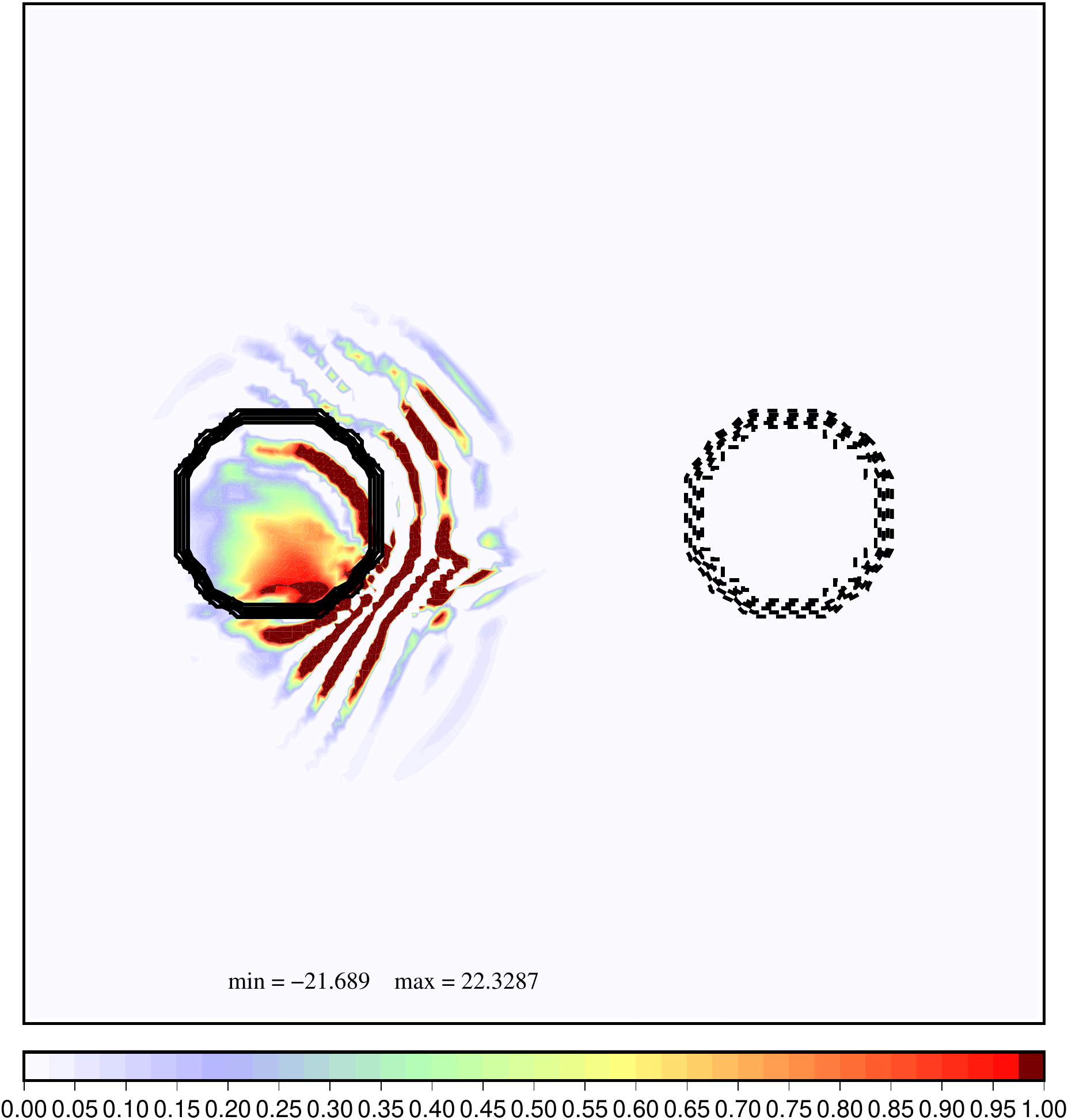} 
  \end{subfigure}
  \begin{subfigure}{0.3\hsize}
    \centering
    {\small (b)}
    \includegraphics[height=52mm, clip, trim=0 40 0 0]{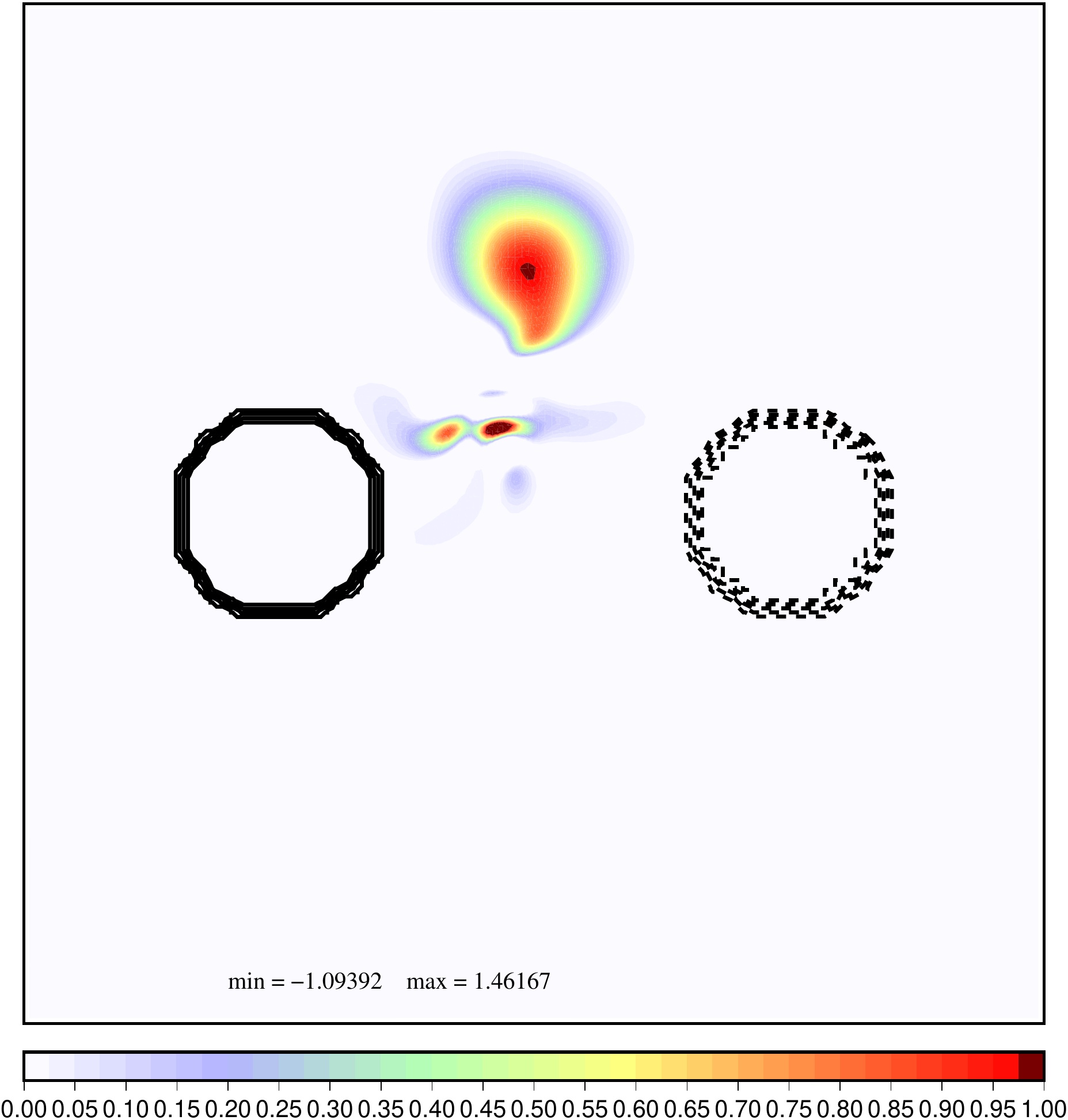} 
  \end{subfigure}
  \begin{subfigure}{0.5\hsize}
    \centering
    \includegraphics[width=70mm, clip, trim=0 0 0 520]{figures/cliff_150.pdf} 
  \end{subfigure}
  \caption{Results of the advection test over steep orography when $A$ is not used in the model. Snapshots are taken at $t$ = (a) 150 s and (b) 600 s. Colour contours show the amplitudes of the tracer density $\rho$. Solid and dashed lines indicate the positive and negative height of orography, respectively, where the contour interval is 100 m.}
  \label{fig:tracer_cliff}
\end{figure}
\begin{figure}[t]
  \centering
  \begin{subfigure}{0.3\hsize}
    \centering
    {\small (a)}
    \includegraphics[height=52mm, clip, trim=0 40 0 0]{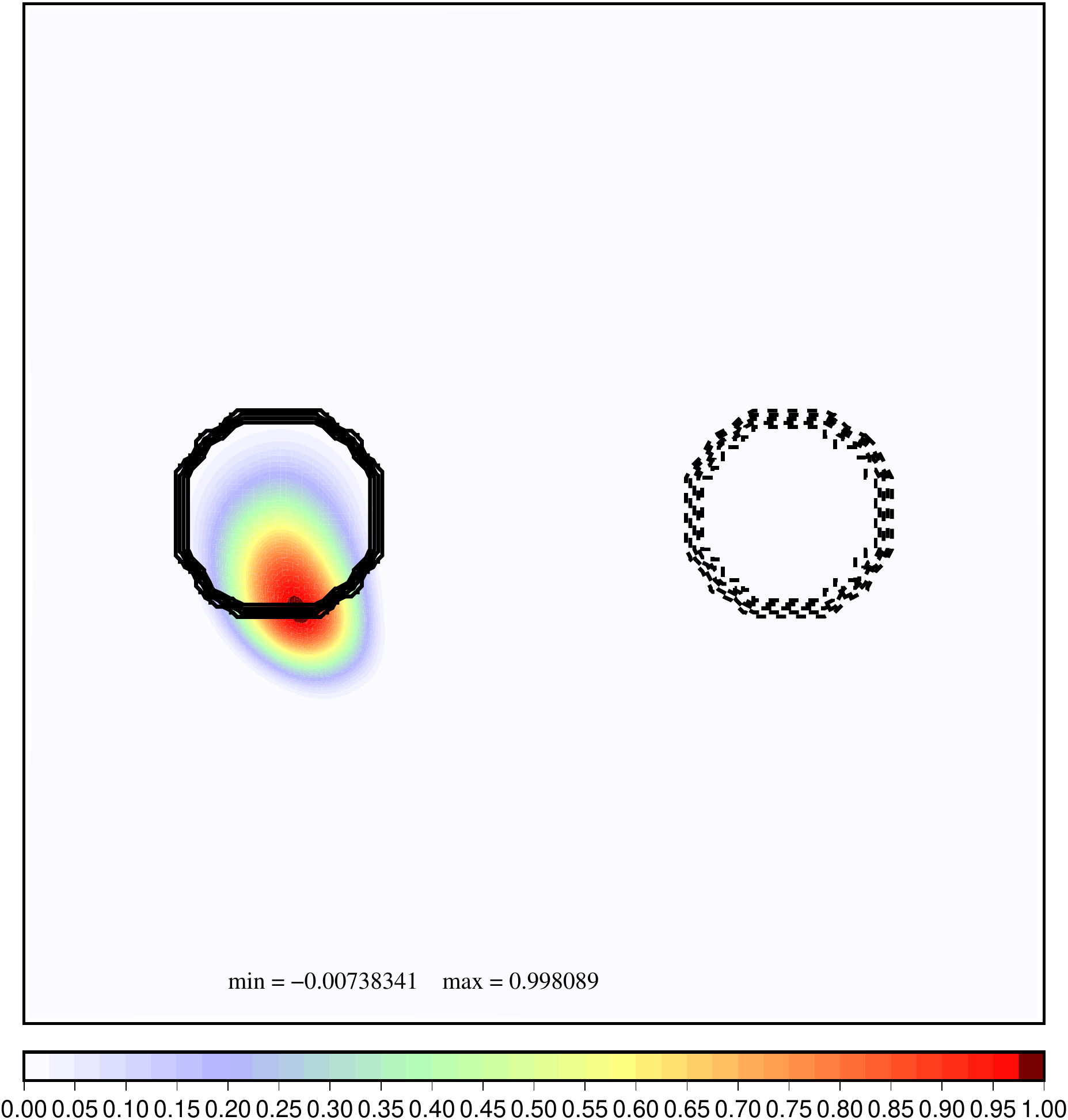} 
  \end{subfigure}
  \begin{subfigure}{0.3\hsize}
    \centering
    {\small (b)}
    \includegraphics[height=52mm, clip, trim=0 40 0 0]{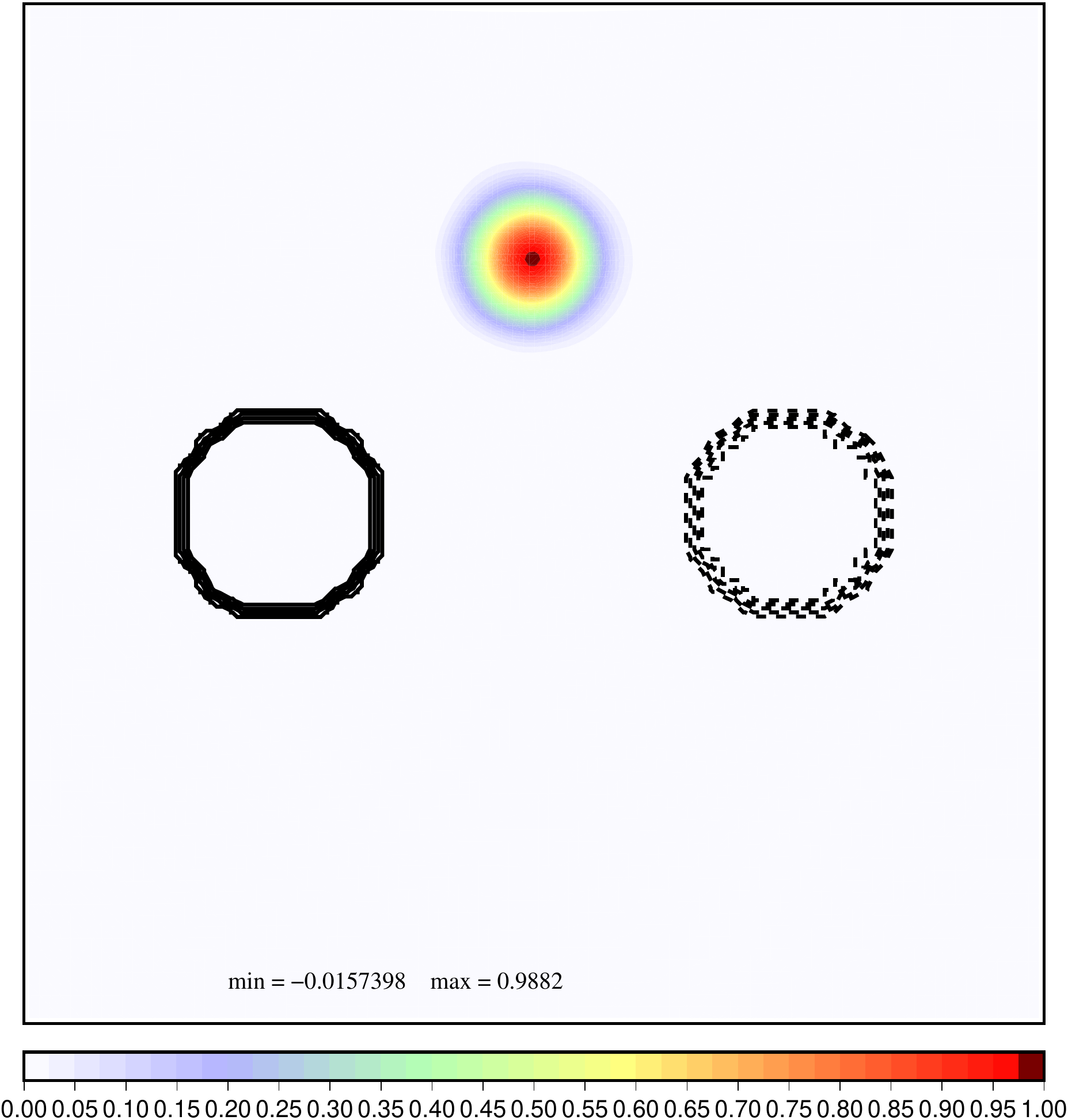} 
  \end{subfigure}
  \begin{subfigure}{0.5\hsize}
    \centering
    \includegraphics[width=70mm, clip, trim=0 0 0 520]{figures/cliffA_150.pdf} 
  \end{subfigure}
  \caption{Same as Figure \ref{fig:tracer_cliff} except that the volume adjustment parameter $A$ is used in the model to adjust the cell volumes.}
  \label{fig:tracer_cliff_A}
\end{figure}

Figure \ref{fig:tracer_cliff} shows the density fields at t = 150 s and t = 600 s when $A$ is not used in the model.
Large numerical errors are found when the tracer passes over the hill (Figure \ref{fig:tracer_cliff}a), and the tracer doesn't recover its original shape after one complete revolution (Figure \ref{fig:tracer_cliff}b).
It might be possible to dampen the oscillations in the solution with an aggressive monotonic advection scheme but monotonic advection should only be achieved for non-divergent wind fields whereas we have introduced artificial divergence by changing the mesh volume. These oscillations should be prevented from occurring by treating the mesh volume correctly.
Figure \ref{fig:tracer_cliff_A} shows the results when $A$ is used to adjust the cell volumes in the model.
In this case, the tracer successfully passes over the hill (Figure \ref{fig:tracer_cliff_A}a) and valley, and completes the revolution without losing shape (Figure \ref{fig:tracer_cliff_A}b). 
Therefore it is shown that the inclusion of $A$ in the model yields stable solutions over steep orography while we observe substantial errors in the model without A, thereby showing the importance of the  maintenance of uniform fields on a moving mesh over orography.

\section{Conclusion}\label{conclusion}

We proposed a novel approach to solve the problem of changes in volume of the domain when resolution changes over orography in a simulation using adaptive meshes.
The volume adjustment parameter is introduced which tracks the true cell volumes by solving an advection equation for a cell volume, achieving conservation of both the volume of the domain and the total mass without tracking the shape of orography within each cell.
The results of a three-dimensional tracer advection test showed that our scheme maintains a uniform field while the mesh resolution changes over orography, whereas the model without the volume adjustment suffers from artificial compression and expansion of the fluid due to the lack of conservation in volume of the domain.
The importance of the maintenance of uniform fields was demonstrated over steep orography where the change in cell volumes on a moving mesh can be pronounced without volume adjustment. The resulting artificial changes in volume lead to large unbounded errors when the advected tracer moves over orography.
The volume adjustment parameter successfully avoided the errors by efficiently tracking the changes in the cell volumes over orography and adjusting the cell volumes. The same idea is considered to be applicable to other variable boundary conditions on a moving mesh (e.g. a land sea mask).
Further work is intended to apply this method to the shallow water or fully compressible equations with the aim of simulating atmospheric problems on a moving mesh over real orography.

\section*{\large Acknowledgement}
We would like to acknowledge NERC grant NE/M013693/1.
The source code for this project is located at \url{https://github.com/AtmosFOAM/AMMM} with the tag \url{paper.2021};
the simulation code itself is located within this repository under \url{run/advection/advectionOTFoam/advectionOT/runAll}.
The codes are developed with OpenFOAM 7 (\url{https://openfoam.org}, 2019).

\begin{appendices}
\section{Boundedness of Volume Adjustments}
\label{appx:bounded}

In section \ref{colin_parameter}, we introduced a volume adjustment parameter, $A$, to correct the cell volumes, $V$, calculated from vertex locations sampled over orography.
As the mesh moves, the sum of all $V$ changes but the sum of $AV$ does not change.
Here we prove that $A$ is bounded above zero which is needed to guarantee that the model may not have negative cell volumes.

The parameter $A$ is calculated from an advection equation \eqref{aeq} discretised using first-order forward in time and first-order downwind in space, which can be rewritten as:
\begin{eqnarray}
\frac{A^{n+1}V^{n+1}-A^{n}V^{n}}{\Delta t} &=&
\sum\limits_{\substack{\rm outward \\ \rm faces}}A_{N}^n \, \underbrace{\phi_{m}}_{\rm positive} \, + \,
\sum\limits_{\substack{\rm inward \\ \rm faces}}A^n \underbrace{\phi_{m}}_{\rm negative} \nonumber \\
&=&
\sum\limits_{\substack{\rm outward \\ \rm faces}}A_{N}^n \, \underbrace{\phi_{m}}_{\rm positive} \, - \,
\sum\limits_{\substack{\rm inward \\ \rm faces}}A^n \underbrace{|\phi_{m}|}_{\rm positive},
\end{eqnarray}
where $A_{N}^n$ denotes the tracer density at the neighbouring cell downstream.
This can be re-arranged for $A^{n+1}$ as a function of values of $A$ at time level $n$:
\begin{eqnarray}
A^{n+1} = \frac{V^{n}}{V^{n+1}}A^{n}\left\{1-\frac{\Delta t}{V^{n}}\sum\limits_{\substack{\rm inward \\ \rm faces}} \underbrace{|\phi_{m}|}_{\rm positive}\right\}
+ \frac{\Delta t}{V^{n+1}} \sum\limits_{\substack{\rm outward \\ \rm faces}}A_{N}^n \, \underbrace{\phi_{m}}_{\rm positive}.
\end{eqnarray}
Given that $A^{n} > 0$ and $A_{N}^n > 0$, we can see that $A^{n+1} > 0$ when
\begin{eqnarray}
1-\frac{\Delta t}{V^{n}}\sum\limits_{\substack{\rm inward \\ \rm faces}} \underbrace{|\phi_{m}|}_{\rm positive} =
1-\underbrace{\left|\frac{\Delta t}{V^{n}}\sum\limits_{\substack{\rm inward \\ \rm faces}} \phi_{m}\right|}_{\rm positive} > 0.
\end{eqnarray}
Since the courant number $C$ is defined as
\begin{eqnarray}
C = \frac{\Delta t}{V^{n}}\sum\limits_{\rm faces}\phi_{m},
\end{eqnarray}
$A^{n+1} > 0$ if $|C| < 1$. Therefore it is proved that, when the initial value of $A$ is positive at all cells, $A$ stays positive as long as the Courant number is less than one.

\section{Numerical Solution of the Monge-Amp\`ere Equation for Mesh Generation}
\label{secn:MAsolution}

The meshing technique is described in full, analysed and compared with other methods full in \cite{browne2016nonlinear} which is summarised here and results are presented for the meshes used in this paper.

\subsection{Introduction}
\label{secn:MAsolution:intro}

An optimally transported mesh is as close as possible to the original mesh (close being defined by the root mean square distance between the vertices of the original and transported mesh) whilst equidistibuting a given scalar monitor function \cite[]{budd2009moving}. To guarantee that the transported mesh is not tangled, the locations, $\bm x$, are defined from the locations of the original mesh, $\bm{\xi}$, by the addition of the gradient of a mesh potential, $\phi$:
\begin{equation}
\bm{x} = \bm{\xi} + \nabla \phi.
\label{eqn:meshPot}
\end{equation}
Equidistribution of the monitor function, $m({\bm x}) > 0$ is expressed as:
\begin{equation}
|\nabla\bm{x}| \ m({\bm x}) = c
\label{eqn:equiDist}
\end{equation}
for a constant $c$ uniform across space where $|\ |$ is the matrix determinant. The combination of equations (\ref{eqn:meshPot}) and (\ref{eqn:equiDist}) gives a fully non-linear elliptic PDE, the Monge-Amp\`ere equation:
\begin{equation}
|I + H(\phi)| = \frac{c}{m(\bm{x})}
\label{eqn:ma3}
\end{equation}
where $I=\nabla\bm{\xi}$ is the identity tensor and $H=\nabla\nabla$ is the Hessian. The meshes in this paper are all the result of numerical solution of the Monge-Amp\'ere equation. 

\cite{budd2009moving} added Laplacian smoothing and a rate of change term to (\ref{eqn:ma3}) making it parabolic and solved using a spectral method. \cite{weller2016mesh} derived an equation to generate optimally transported meshes on the surface of a sphere, linearised about a uniform flat mesh to create fixed point iterations, each iteration requiring the solution of a Poisson equation discretised using finite volumes. \cite{McRae2018} re-wrote the equation on the surface of a sphere as a PDE and solved using a Newton solver with finite elements. Here we describe a Newton method for solving the Monge-Amp\`ere equation on a finite plane and discretise in space with finite volumes following \cite{weller2016mesh}.

\subsection{Numerical Method}
\label{secn:MAsolution:method}

We define a Newton method for solving (\ref{eqn:ma3}) in Euclidean geometry, linearising the LHS around the previous iteration and using the RHS from the previous iteration. $\bm {x}^k = \bm{\xi} + \nabla \phi^k$ is the solution at iteration $k$. By writing $\phi^{k+1} = \phi^k+\varepsilon\psi$ it can be shown that
\begin{equation}\label{delaware}
|I+H(\phi^{k+1})| = |I+H(\phi^k)| + \nabla \cdot P^k\nabla \varepsilon \psi + \mathcal{N}(\varepsilon \psi),
\end{equation}
where $P^k$ is the matrix of cofactors of $I+H(\phi)$ and $\mathcal{N}$ is some nonlinear function. In 2D
\begin{equation}
P^k = \begin{bmatrix} 1 + \phi^k_{yy} & -\phi^k_{xy} \\ -\phi^k_{xy} & 1+\phi^k_{xx}\end{bmatrix}
\end{equation}
and $\mathcal{N}(\varepsilon \psi) =  \varepsilon^2 |H(\psi)|$.
In 3D, a more involved calculation can show $\mathcal{N}(\varepsilon \psi) =  \varepsilon^3 \mathcal{\tilde{N}}(\psi)$ and 
\begin{equation*}
P^k = \begin{bmatrix}
1 + \phi^k_{yy} + \phi^k_{zz} + \phi^k_{yy}\phi^k_{zz} - \phi^k_{yz}\phi^k_{yz} & -\phi^k_{xy} -\phi^k_{xy}\phi^k_{zz} + \phi^k_{xz}\phi^k_{yz} & -\phi^k_{xz} -\phi^k_{xz}\phi^k_{yy}+\phi^k_{xy}\phi^k_{yz}\\
-\phi^k_{xy} -\phi^k_{xy}\phi^k_{zz} +\phi^k_{xz}\phi^k_{yz} & 1+\phi^k_{xx} + \phi^k_{zz} + \phi^k_{xx}\phi^k_{zz} - \phi^k_{xz}\phi^k_{xz}&-\phi^k_{yz} - \phi^k_{xx}\phi^k_{yz}+\phi^k_{xy}\phi^k_{xz}\\
-\phi^k_{xz} - \phi^k_{xz}\phi^k_{yy} + \phi^k_{xy}\phi^k_{yz} & -\phi^k_{yz} - \phi^k_{xx}\phi^k_{yz} + \phi^k_{xy}\phi^k_{xz} & 1+ \phi^k_{xx}+\phi^k_{yy}+\phi^k_{xx}\phi^k_{yy}-\phi^k_{xy}\phi^k_{xy}
\end{bmatrix}.
\end{equation*}
At convergence terms proportional to $\varepsilon^d$ (where $d$ is the dimensionality of space) will disappear so at each iteration, we solve the following Poisson equation for $\varepsilon\psi$:
\begin{equation}
\nabla \cdot \left(P^k\nabla \varepsilon \psi\right) = \frac{c}{m(\bm{x}^k)}
 - |I+H(\phi^{k})| .
\label{eqn:MA_k}
\end{equation}
Equation (\ref{eqn:MA_k}) is elliptic as long as $P^k$ is positive definite. For simplicity and efficiency we use a finite volume discretisation for spatial discretisation of (\ref{eqn:MA_k}). However, unlike wide stencil finite difference methods \cite[e.g.][]{oberman2008wide}, this is not guaranteed to give monotonic solutions . This means that $P^k$ can become non-positive definite so (\ref{eqn:MA_k}) loses its ellipticity and solutions rapidly diverge. To remedy this we can modify \eqref{eqn:MA_k} to maintain ellipticity by replacing the matrix $P^k$ with a modified matrix $Q^k$ such that
\begin{equation}\label{eqn:fp3gamma}
Q^k = P^k + \gamma I
\end{equation}
and $\gamma$ is defined as
\begin{equation}\label{eqn:eigenvalue}
\gamma := \begin{cases}
0 \qquad &\text{if} \qquad \min\sigma[P^k] > 0\\
\delta - \min \sigma[P^k] \qquad &\text{if} \qquad \min\sigma[P^k] \le 0.
\end{cases}
\end{equation}
The constant $\delta>0$ is chosen to avoid round-off errors (we have taken $\delta=10^{-5}$), and $\sigma[P^k]$ refers to the spectrum of $P^k$. This process simply shifts the eigenvalues of the matrix $P^k$ so that they remain positive.

The iterations labelled $k$ are called outer iterations because the Poisson equation is also solved using an iterative solver within each outer iteration. 

The Laplacian and the Hessian of (\ref{eqn:MA_k}) are discretised in space using compact finite volumes, following \cite{weller2016mesh}. This is equivalent to second order finite differences on a uniform grid. Zero gradient boundary conditions are used. The spatial discretisation leads to a set of linear simultaneous equations. These are solved using the OpenFOAM GAMG solver with a symmetric Gauss Seidel smoother and an LU pre-conditioner. A maximum of 10 solver iterations are allowed. The tightest solver tolerance is $10^{-4}$ but the solver is only solved to a tolerance of 0.01 times the initial residual each outer iteration. This is to avoid spending too much time solving the first few iterations tightly when subsequent iterations will have updated coefficients.

\subsection{The Monitor Function}
\label{secn:MAsolution:monitor}

The monitor function is based on the Frobenius norm of the Hessian of the tracer density, $\rho$, which in two dimensions is
\begin{equation}
m_1(\bm{x}) = \sqrt{\rho_{xx}^2 + \rho_{xy}^2 + \rho_{yx}^2 + \rho_{yy}^2}.
\label{eqn:monitorGradGrad}
\end{equation}
Following \cite{McRae2018} we use the rule of thumb that half of the resolution should be placed where not much is happening. This can be approximately achieved by setting:
\begin{equation}
m_2(\bm{x}) = \min\biggl(\frac{m_1}{\overline{m_1}+1}, r_{\max}\biggr)
\label{eqn:monitor5050}
\end{equation}
where $\overline{m_1}$ is the area average of $m_1$ and $r_{\max}$ is the ratio of smallest to largest cell volumes/areas of the resulting adapted mesh. For the simulations in section \ref{results} we use $r_{\max}=4$ meaning that, if cells have aspect ratio 1 then the maximum ratio of smallest to largest cell side lengths is 2. The monitor function is smoothed before it is used for mesh generation so that the resulting mesh varies smoothly, which is advantageous for finite volume and finite difference methods that have the property of super convergence. The final monitor function, $m_3$, is the implicit solution of the diffusion equation:
\begin{equation}
\frac{m_3 - m_2}{\Delta t} = K \nabla^2 m_3
\label{eqn:monitorSmooth}
\end{equation}
where the diffusion coefficient, $K$, is mesh size and time step dependent:
\begin{equation}
K = M\frac{\Delta x^2}{4\Delta t}.
\label{eqn:monitor3}
\end{equation}
$M$ is equivalent to the number of applications of a $(1,-2,1)$ filter to smooth the monitor function. $M=20$ is used for the meshes presented in section \ref{results}. The Laplacian in (\ref{eqn:monitorSmooth}) is calculated on the uniform orthogonal computational mesh of squares using 2nd-order centred differences (i.e. $(1,-2,1)$ differencing in each direction).

\subsection{Results}
\label{secn:MAsolution:results}

Meshes are generated for the linear advection results using a cosine-shaped tracer in section \ref{doubleCosine} starting from initial uniform grids of $50\times 50$, $100\times 100$, $200\times 200$ and $400\times 400$ points in a plane of size 10 km by 10 km.

\begin{figure}[h]
\begin{tabular}{cc}
\multicolumn{2}{c}{Convergence of mesh generation starting from regular mesh}\vspace{0.5em}\\
\includegraphics[width=0.48\linewidth]{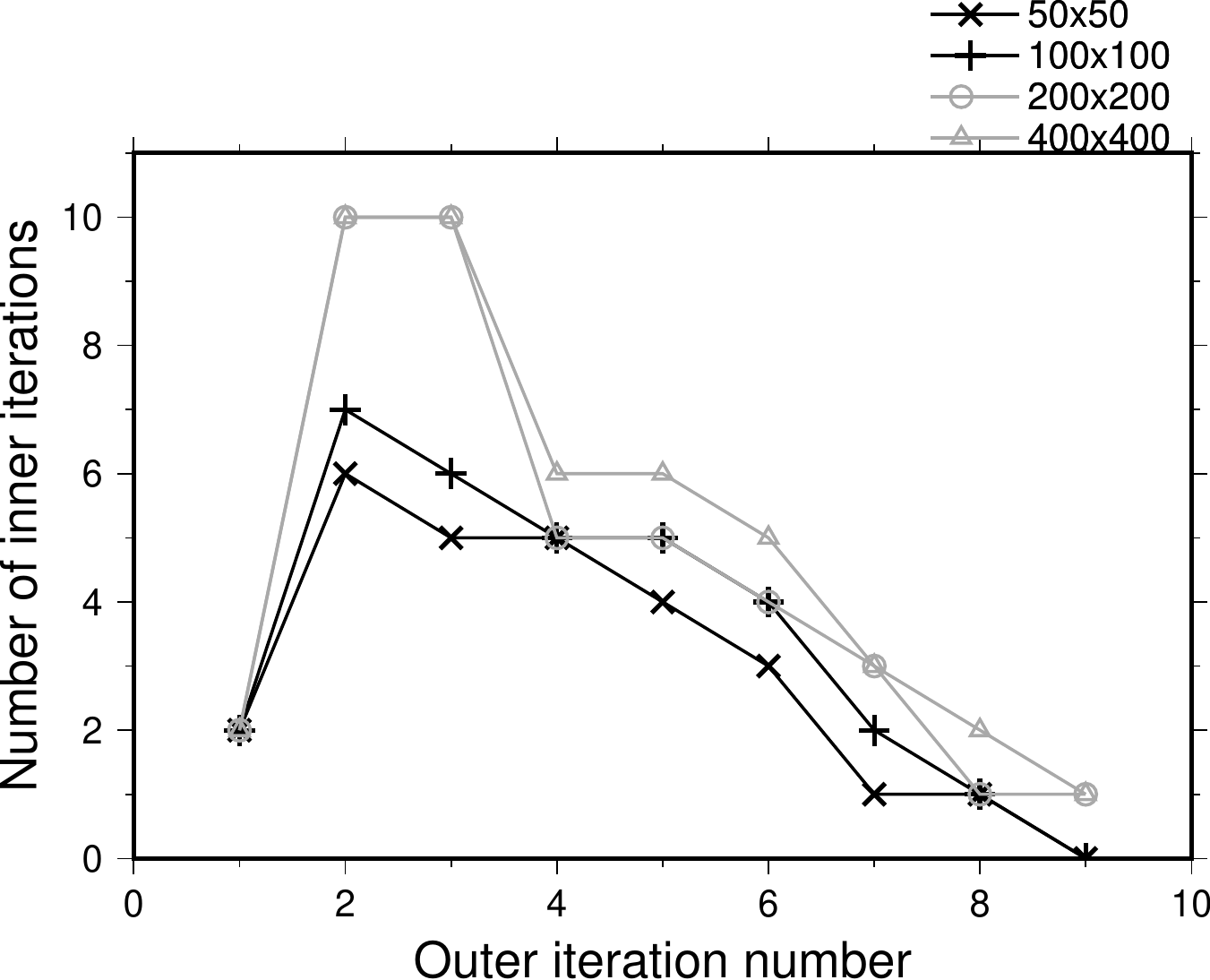}
&
\includegraphics[width=0.48\linewidth]{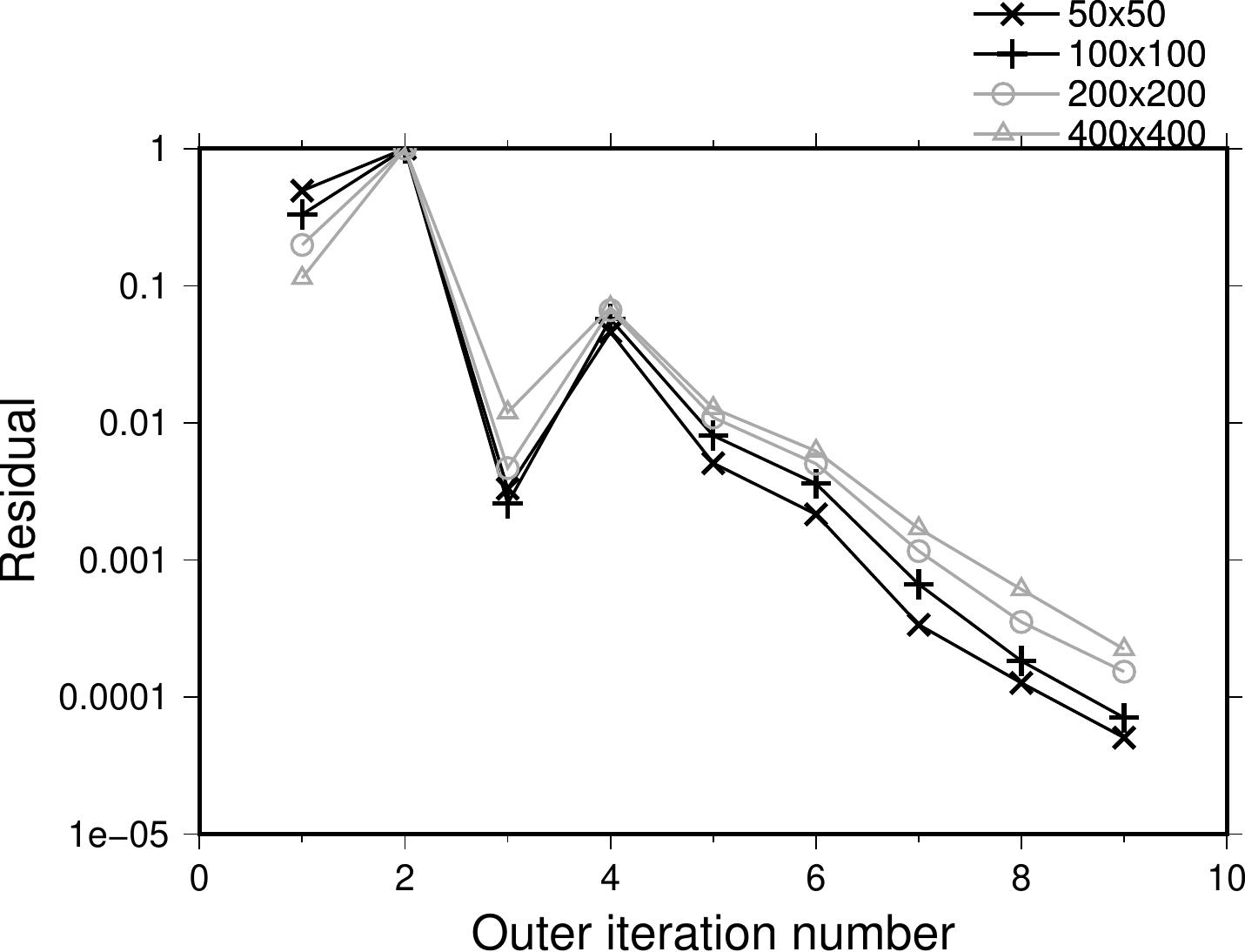}\vspace{0.5em}
\\
\multicolumn{2}{c}{Convergence of mesh generation each time step}\\
\includegraphics[width=0.48\linewidth]{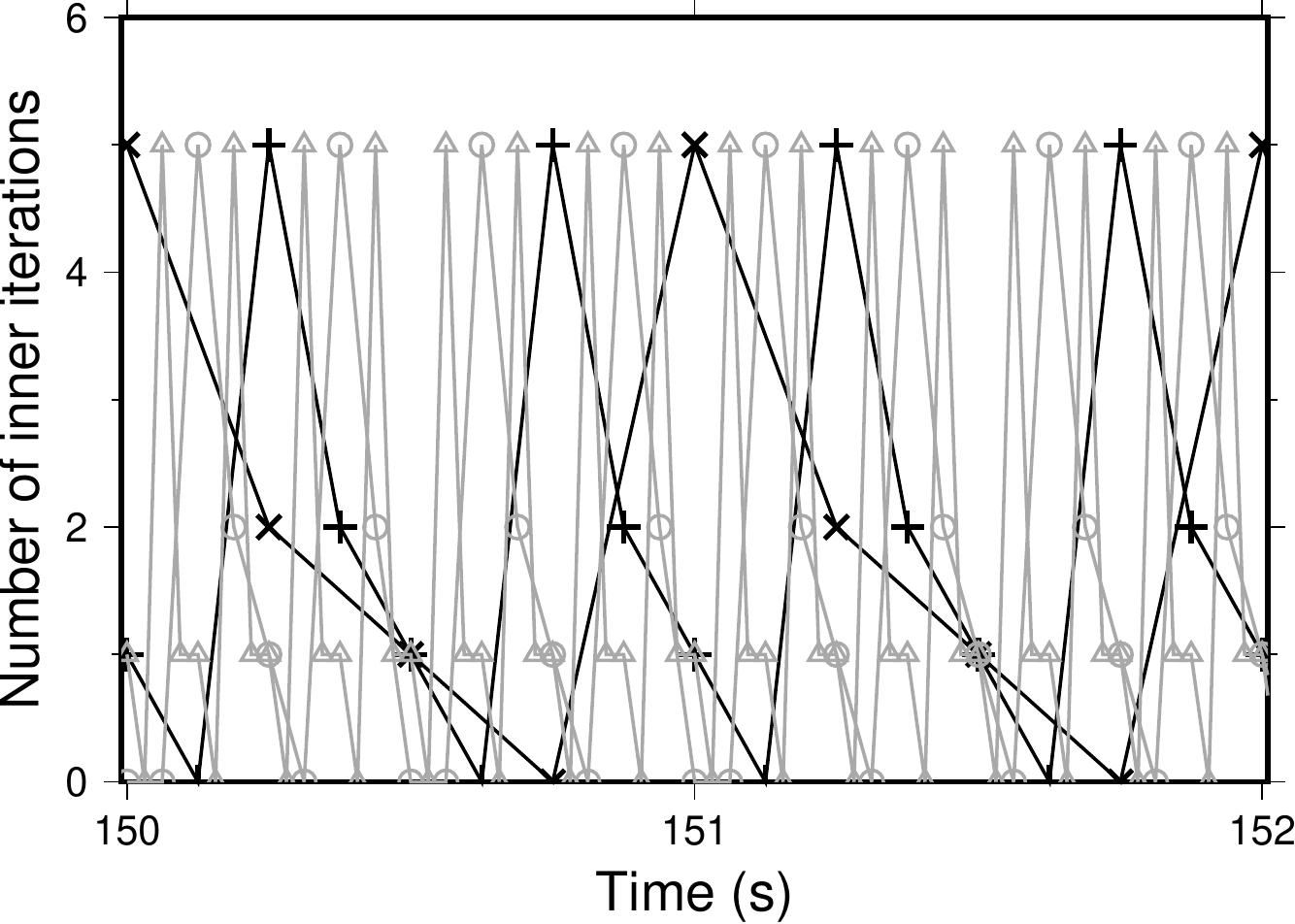}
&
\includegraphics[width=0.48\linewidth]{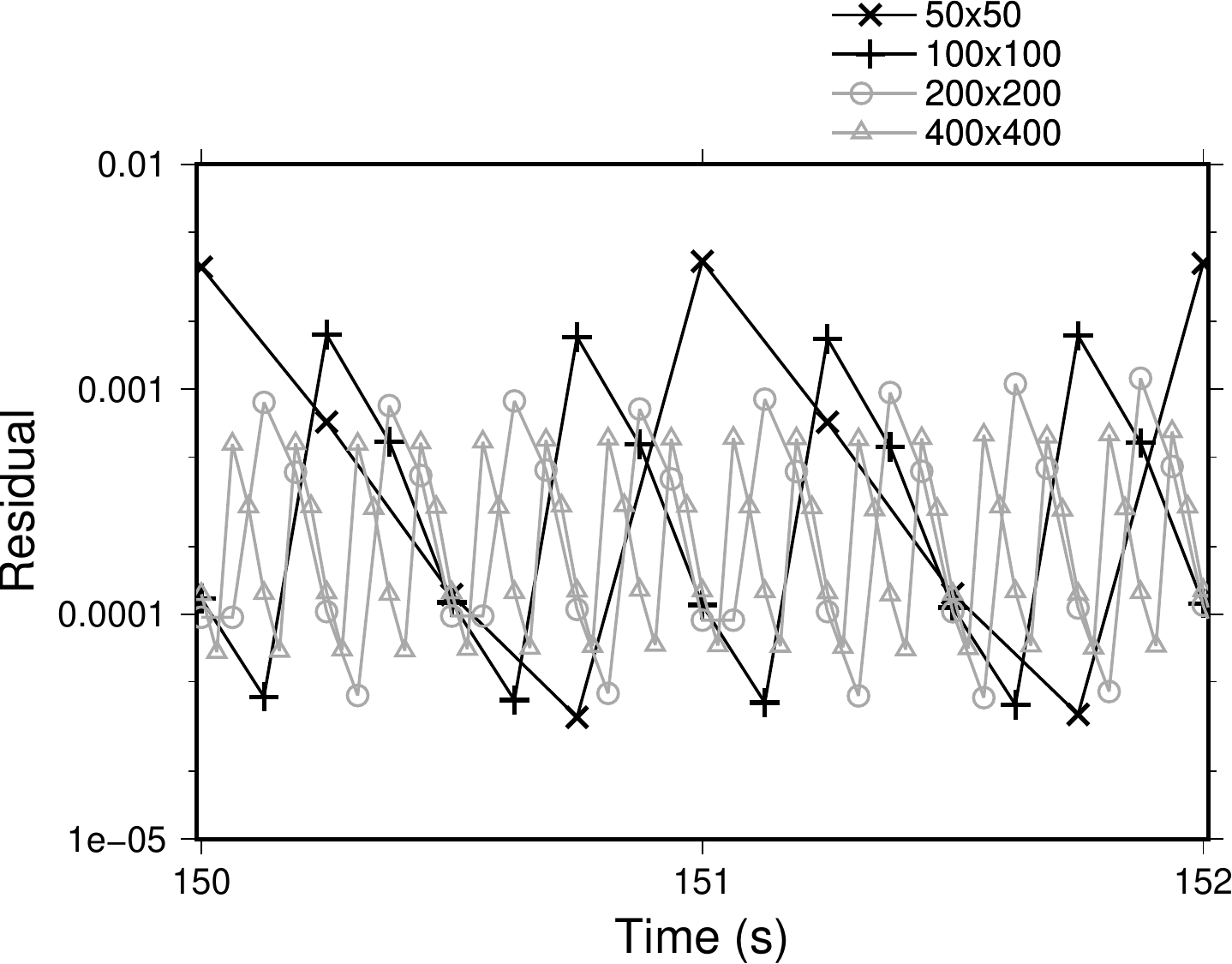}
\\
\end{tabular}
\caption{Convergence of mesh generation for four different resolutions for advection over a hill and a valley. The top row shows convergence starting from a regular mesh. The bottom row shows convergence for each time step of the transient simulation. In the transient simulation there are four outer iterations, each consisting of one solution of the Poisson equation (\ref{eqn:MA_k}) initialised from the mesh at the previous time step. The left hand side shows the number of linear equation solvers per solution of the Poisson equation and the right hand side shows the initial residual before the linear equation solver is called. Before the simulation is started, a refined mesh is calculated from a uniform mesh and convergence is shown in the top row.}
\label{fig:meshConverge}
\end{figure}

Before the advection simulation starts, an initial mesh is generated using the monitor function calculated from the analytic description of the initial conditions. 9 outer iterations are used. The residual of the Poisson equation solver and the number of iterations of the Poisson equation solver for each outer iteration are shown in the top row of Figure \ref{fig:meshConverge} for all resolutions. Convergence is reasonably insensitive to resolution which is necessary for efficiency. Convergence in the first three iterations is noisy but then convergence proceeds exponentially (note the residuals are on a log-scale).

While solving the advection equation, the mesh is moved every time step. The same uniform, regular computational mesh is used to solve the Monge-Amp\`ere equation each time step but the solution is initialised from the previous time step. Each time step, a maximum of 4 outer iterations of the Monge-Amp\`ere Newton solver are allowed. The number of inner (linear equation solver) iterations and the initial residual for each solver are shown in the bottom row of Figure \ref{fig:meshConverge}. This shows that at most 5 inner iterations are needed and convergence is exponential between each outer-iteration per time step.

\subsection{Further Remarks}
\label{secn:MAsolution:remarks}

Section \ref{secn:MAsolution:method} described a regularisation technique to ensure that the discretised Poisson equation (\ref{eqn:MA_k}) remains elliptic by artificially increasing the diagonal of the Poisson equation coefficient, $P$. This may raise concerns that we are arbitrarily changing the problem that we are solving. However the regularisation is only very occasionally needed and is only ever needed during the first one or two outer iterations and so this regularisation never influences the final converged solution.

\cite{browne2016nonlinear} compared this Newton solver with the parabolic method of \cite{Browne2014a} and with the fixed point iterations used by \cite{weller2016mesh}. Convergence of the proposed Newton solver was far superior and free of arbitrary parameters. \cite{browne2016nonlinear} also proposed a Newton solver that involved linearising the $c/m$ term of (\ref{eqn:MA_k}). This lead to even faster but unreliable convergence and so is not used here.

\end{appendices}

\bibliography{advectionMM}

\end{document}